\documentclass[12pt]{article}

\usepackage{fixmath}
\usepackage{amsmath}
\usepackage{amsfonts}
\usepackage{amssymb}
\usepackage{dsfont}
\usepackage{eucal}
\usepackage{mathrsfs}
\usepackage{verbatim}
\usepackage{mathtools}

\newcommand{\tang}{\mathbf{T}}
\newcommand{\hu}{\hat{u}}
\newcommand{\wD}{{\sf D}}
\newcommand{\hF}{{\hat{F}}}
\newcommand{\hphi}{{\hat{\varphi}}}
\newcommand{\hnu}{{\hat{\nu}}}
\newcommand{\hgamma}{{\hat{\gamma}}}
\newcommand{\htau}{{\hat{\tau}}}
\newcommand{\sbt}{\,\begin{picture}(-1,1)(-1,-3)\circle*{2}\end{picture}\ }
\newcommand{\dder}{{\sbt}}
\newcommand{\hf}{\frac{1}{2}}
\newcommand{\ric}{{\operatorname{Ric}}}
\newcommand{\R}{\mathbb{R}}  
\newcommand{\RE}{\mathbb{R}}  
\newcommand{\Bx}{\mathbf{x}}  
\newcommand{\By}{\mathbf{y}}  
\newcommand{\Bdelta}{\mathbold{\delta}}  
\newcommand{\F}{\mathcal{F}}  
\newcommand{\N}{\mathbb{N}}  
\newcommand{\Z}{\mathbb{Z}}  
\newcommand{\D}{\mathbb{D}}  
\newcommand{\I}{\mathds{1}}  
\newcommand{\EF}{\mathcal{E}} 
\newcommand{\eps}{\varepsilon} 
\newcommand{\dSi}{\partial_{S_i}} 
\newcommand{\CM}{\mathbb{H}}  
\newcommand{\HE}{\mathcal{H}}  
\newcommand{\cyl}{\mathcal{F}C_b^\infty}
\newcommand{\nti}{{n\to \infty}}
\newcommand{\lni}{\underset{n\rightarrow \infty}{\longrightarrow}}
\newcommand{\lki}{\underset{k\rightarrow \infty}{\longrightarrow}}

 \usepackage{amsthm}

\numberwithin{equation}{section}

\theoremstyle{definition}

\newtheorem{thm}{Theorem}[section]

\newtheorem{proposition}[thm]{Proposition}

\newtheorem{lemma}[thm]{Lemma}

\newtheorem{corollary}[thm]{Corollary}

\newtheorem{remark}[thm]{Remark}

\begin{document}
\title{A class of infinite dimensional stochastic processes with unbounded diffusion} 
\author{John Karlsson and J\"org-Uwe L\"obus\\
 \normalsize Matematiska institutionen \\ [-1.2ex]
 \normalsize Link\"opings universitet \\ [-1.2ex]
 \normalsize SE-581 83 Link\"oping \\ [-1.2ex]
 \normalsize Sverige
}
\date{}
\maketitle
{\footnotesize
\noindent
\begin{quote}
\oddsidemargin5.0cm
\textwidth10.0cm

{\bf Abstract}
The paper studies Dirichlet forms on the classical Wiener space and the Wiener space over non-compact complete Riemannian manifolds. 
The diffusion operator is almost everywhere an unbounded operator on the Cameron--Martin space. 
In particular, it is shown that under a class of changes of the reference measure, quasi-regularity of the form is preserved. We also show that under these changes of the reference measure, derivative and divergence are closable with certain closable inverses. We first treat the case of the classical Wiener space and then we transfer the results to the Wiener space over a Riemannian manifold.

\noindent


\end{quote}
}

\section{Introduction}

This paper is concerned with Dirichlet forms of type
\begin{align}
\begin{split}
\label{first}
\EF (F,G) &= \int \left\langle DF,ADG \right\rangle_\CM \, \varphi
d\nu \\
&=\int \left\langle DF,\sum_{i=1}^\infty \lambda_i \left\langle S_i,DG\right\rangle_\CM S_i \right\rangle_{\CM}\varphi d\nu,
\end{split}
\end{align}
where the diffusion operator $A$ is in general $\nu$-a.e. unbounded, cf. \cite{Lobus2004}. This topic has received increasing interest over the past years, see e.g. \cite{ChenWu2014,WangWu2008,WangWu2009}. We are interested in weight functions $\varphi$ of the form
\begin{equation*}
\varphi(\gamma)=\exp \left\{ \int _0^1 \langle b(\gamma_s), \, d\gamma_s \rangle_{\RE^d} - \hf \int_0^1 |b(\gamma_s)|^2 \, ds \right\}.
\end{equation*}
This choice of weight functions is motivated by the geometric setting in the papers \cite{WangWu2008,WangWu2009} by F.-Y. Wang and B. Wu. Accordingly the weight function $\varphi$ will be specified in Section \ref{section_geometric}, where we consider the above bilinear form in a geometric framework. In addition we study a weight function of similar type in the flat case. We present our ideas in terms of infinite dimensional processes on the classical Wiener space using a coordinate representation. This makes the subject comprehensible, in particular, to readers familiar with the finite dimensional theory. The form is studied on the set of smooth cylindrical functions of type
\begin{equation*}
F,G\in Y=\{F\left(\gamma)=f(\gamma(s_1),\dots,\gamma(s_k)\right):s_j \text{ is a dyadic point}\},
\end{equation*}
where $\gamma$ is a Wiener trajectory. Simultaneously we will also study
the form on the more common set of cylindrical functions
\begin{equation*}
F,G\in Z=\{F\left(\gamma)=f(\gamma(s_1),\dots,\gamma(s_k)\right):s_j \in [0,1]\}.
\end{equation*}
Well-definiteness of $\EF$ on $Y$ is a consequence of the fact that
the sum in (\ref{first}) is finite. Well-definiteness of $\EF$ on $Z$
requires the convergence of the sum in (\ref{first}). These two
different initial situations result in possibly different closures of
$(\EF,Y)$ and $(\EF,Z)$ on $L^2(\varphi \nu)$. Using the coordinate
representation in (\ref{first}) we give conditions for closability. It
turns out that in the classical case i.e. $\lambda_1=\lambda_2=\ldots=1$, the condition $\varphi^{-1} \in L^1(\nu)$ is
sufficient for closability. However, the paper investigates forms of
the structure (\ref{first}) with an in general $\nu$-a.e. unbounded
diffusion operator. We give necessary and sufficient conditions  on
$\lambda_1(\gamma), \lambda_2(\gamma), \dots$ that guarantee
closability of $(\EF,Y)$ and $(\EF,Z)$ on $L^2(\varphi \nu)$ in terms
of the Schauder functions $S_i$, $i\in \N$, the coordinate functions
 in the Cameron--Martin space. Locality, Dirichlet property, and quasi-regularity of the closure $(\EF,Z)$ on $L^2(\varphi \nu)$ is then obtained by using methods of \cite{BouleauHirsch,Driver1992,Lobus2004,Rockner,Schmuland1992}.

In the paper \cite{Driver1992} the classical quasi-regular Dirichlet form
on a compact Riemannian manifold has been constructed and associated with
a diffusion process.
\\The paper \cite{Lobus2004} extends the same idea in particular to the case of unbounded diffusion. There the form is
\begin{equation}
\label{eq_form_study}
\int \langle DF,ADG \rangle_\CM \, d\nu
\end{equation}
where the operator $A$ acts in $L^2(\nu;\CM)$, has a representation
\begin{equation*}
A\Phi(\gamma):=\sum_{i=1}^\infty \lambda_i \langle S_i,\Phi(\gamma)\rangle_\CM S_i,
\end{equation*}
and the constant diffusion coefficients $\lambda_n$ are assumed to be bounded from below and possibly unbounded. It is shown that under the assumption $\lambda_n \le cn^{1-\varepsilon}$, the form is a quasi-regular Dirichlet form and thus has an associated process.\\
\\The approach in the present paper is similar to that in
\cite{WangWu2008}. Again, infinite dimensional diffusion
processes are investigated, where the diffusion is related to possibly unbounded
operators. The present paper and \cite{WangWu2008} refer to the ideas of
\cite{Lobus2004}, however study forms of type \eqref{eq_form_study} with
non-constant diffusion coefficients. The theory in \cite{WangWu2008} is
presented in a geometric framework, i.e. on the Wiener trajectories
over a complete non-compact Riemannian manifold while here we first consider the form on the classical Wiener space and
later turn to a geometric setting. In addition we also take
particular consideration to the special case when the form is defined
pointwise, constant over the trajectory space as in \cite{Lobus2004}. \\
The operator $(A,D(A))$ in \cite{WangWu2008} is introduced as a densely defined self-adjoint operator such that
\begin{align}
\tag{A0}
\begin{split}
&Af\Phi=fA\Phi \text{ for } f\in L^\infty(\nu) \text{ and } \Phi\in D(A)\text{ such that }f\Phi\in D(A) \\
&\text{and } A\ge \eps\operatorname{Id} \text{ for some }\varepsilon>0.
\end{split}
\end{align}
Here $\nu$ is the Wiener measure on the flat trajectory space,
or trajectory space over a complete non-compact Riemannian manifold. Furthermore
$\operatorname{Id}$ denotes the identity operator. Under (A0) the following is obtained. For $f,\Phi$ as above and $\Psi\in D(A)$ such that
$f\Psi\in D(A)$ it holds that
\begin{align*}
  \int f\int_0^1 (A\Phi)^\dder \cdot \dot{\Psi}\, ds\,
  d\nu&=\int\int_0^1\dot{\Phi}(A(f\Psi))^\dder\, ds\, d\nu\\
&=\int f\int_0^1
  \dot{\Phi}(A\Psi)^\dder\, ds\,d\nu,
\end{align*}
which implies
\begin{equation*}
  \langle A\Phi(\gamma),\Psi(\gamma)\rangle_\CM=\langle
  \Phi(\gamma),A\Psi(\gamma)\rangle_\CM,
\end{equation*}
for $\nu$-a.e. $\gamma\in\Omega$. This observation suggests that
\begin{equation*}
(A\Phi)(\gamma)=A(\gamma)\Phi(\gamma),
\end{equation*}
cf. assumption (2) in Remark 1.1 of \cite{WangWu2008}, that is $(A\Phi)^{1/2}(\gamma)=A(\gamma)^{1/2}\Phi(\gamma)$ for
$\nu$-a.e. $\gamma\in\Omega$. Below we specify
\begin{equation*}
D(A):=\left\{\Phi \in L^2(\nu;\CM): \int \sum_{i=1}^\infty \lambda_i(\gamma)^2  \left\langle S_i,\Phi(\gamma) \right\rangle_{\CM}^2 \, d\nu <\infty\right\},
\end{equation*}
\begin{align*}
A\Phi(\gamma):=\sum_{i=1}^\infty \lambda_i(\gamma)  \left\langle S_i,\Phi(\gamma) \right\rangle_{\CM}S_i ,\quad \gamma \in \Omega, \quad \Phi \in D(A).
\end{align*}

In \cite{WangWu2008} the form
\begin{equation*}
\EF_A(F,G)=\int_{W_0} \langle A^{1/2}DF,A^{1/2}DG\rangle_\CM \, d\mu, 
\end{equation*}
is defined where $F,G\in \{f\in \cyl :Df\in D(A^{1/2})\}$ and $\cyl$ denotes the set of cylindrical functions
$\{F\left(\gamma)=f(\gamma(s_1),\dots,\gamma(s_k)\right):s_j \in
[0,1], f\in C_b^\infty\}$. We take over this concept however in the flat case we use
the cylindrical functions in $Y$ or $Z$ for the initial definition. 

We give a necessary and sufficient condition for well-definiteness and closability of the bilinear form. In the particular case of non-decreasing eigenvalues  $0<\lambda_1(\gamma)<\lambda_2(\gamma),\ldots$ it reads as 
	\begin{equation*}
		\sum_{p=0}^\infty  \frac{\int \lambda_{d2^p} \, \varphi d\nu}{2^p} < \infty.
	\end{equation*}

In order to handle weights, $\varphi$, with respect to the Wiener measure, the paper \cite{WangWu2008} introduces geometric conditions. We use
\begin{equation*}
	\frac{f}{\varphi} \in L^1(\nu),\quad f\in Z,
\end{equation*}
to verify closability of the classical form, which was first treated in \cite{Driver1992}. To show quasi-regularity this paper and \cite{WangWu2008} both use the method of \cite{Driver1992}, \cite{Lobus2004}, and \cite{Schmuland1992}. The proof in \cite{WangWu2008} involves a geometric result of \cite{Wang2004}. We use alternative conditions on the manifold.

\section{Definitions}

We study the form on the space $L^2(\varphi \nu)\equiv L^2(\Omega, \varphi \nu)$ where $\Omega:=C_0([0,1];\RE^d):=\{\gamma\in C([0,1]; \RE^d, \gamma(0)=0\}$, $\nu$ is the Wiener measure on $\Omega$, and $\varphi$ is a density function specified below. As stated earlier, the form is given by
\begin{equation}
\EF (F,G) = \int \left\langle DF,ADG \right\rangle_\CM \, \varphi d\nu, \quad F,G\in D(\EF),
\end{equation}
where $\CM$ is the Cameron-Martin space, i.e. the space of all absolutely continuous $\RE^d$-valued functions $f$ on $[0,1]$, with $f(0)=0$ and equipped with inner product
\begin{equation*}
\langle \varphi,\psi\rangle_\CM:=\int_{[0,1]} \langle  \varphi'(x) , \psi'(x) \rangle_{\RE^d} \,dx.
\end{equation*}
Motivated by \cite{WangWu2008}, we suppose for proving quasi-regularity of the form that  $\varphi~:~\Omega\rightarrow [0,\infty]$ has the form
\begin{equation}
\label{def_phi}
\varphi(\gamma)=\exp \left\{ \int _0^1 \langle b(\gamma_s), \, d\gamma_s \rangle_{\RE^d} - \hf \int_0^1 |b(\gamma_s)|^2 \, ds \right\}
\end{equation} 
where we choose $b$ from two different points of view. On the one hand we are interested in closability, see Section \ref{section_closability}. On the other hand in Section \ref{section_geometric}, these coefficients correspond to a change of the Wiener measure to a certain diffusion measure on the manifold. We define the set of all cylindrical functions
\begin{align*}
Z:=&\Big\{F(\gamma)=f\left(\gamma(s_1),\dots,\gamma(s_k)\right), \gamma \in \Omega : \\
& 0<s_1<\dots<s_k=1, f\in C^{\infty}_p\left(({\RE^d})^k\right), k\in \N\Big\}
\end{align*}
where $f\in C^{\infty}_p$ means that $f$ and all its partial derivatives are smooth with polynomial growth. We also define
\begin{align*}
Y:=&\Big\{F(\gamma)=f\left(\gamma(s_1),\dots,\gamma(s_k)\right), \gamma \in \Omega : \\
&F \in Z, s_1,\dots,s_k\in \left\{\textstyle{\frac{l}{2^n}}:l \in \{1,\dots,2^n\}\right\}, n\in \N \Big\}.
\end{align*}
For $F\in Z$ and $\gamma \in \Omega$ the gradient operator $D$ is defined by 
\begin{equation}
\label{def_grad}
D_sF(\gamma)=\sum_{i=1}^{k} (s_{i}\wedge s) (\nabla_{s_{i}}f)(\gamma), \quad s\in [0,1],
\end{equation}
where $(\nabla_{s_{i}}f)(\gamma)=(\nabla_{s_{i}}f)(\gamma(s_1),\dots,\gamma(s_k))$ denotes the gradient of the function $f$ relative to the $i$th variable while holding the other variables fixed. 
We let $(e_j)_{j=1,\dots,d}$ denote the standard basis in $\RE^d$ and
\begin{align}
\begin{split}
\label{eq_def_haar1}
&H_1(t)=1, \\ 
  &H_{2^m+k}(t) = \left\{
  \begin{array}{l l l}
    2^{m/2} & \quad \text{if $t \in \left[\frac{k-1}{2^m}, \frac{2k-1}{2^{m+1}}\right)$}\\
    -2^{m/2} & \quad \text{if $t \in \left[\frac{2k-1}{2^{m+1}}, \frac{k}{2^m}\right)$} \qquad k=1,\dots,2^m, \, m=0,1,\dots,\\
    0 & \quad \text{otherwise}
  \end{array} \right.
  \end{split}
\end{align}
denote the system of the Haar functions on $[0,1]$. We also define 
\begin{equation}
\label{eq_def_gn}
g_{d(r-1)+j}:=H_r \cdot e_j, \quad r\in \N, \, j\in \{1,\dots,d\},
\end{equation}
and
\begin{equation*}
S_n(s):=\int_0^s g_n(u) \, du, \quad s\in [0,1],\, n\in\N.
\end{equation*}

Let
\begin{equation*}
D(A):=\left\{\Phi \in L^2(\nu;\CM): \int \sum_{i=1}^\infty \lambda_i(\gamma)^2  \left\langle S_i,\Phi(\gamma) \right\rangle_{\CM}^2 \, d\nu <\infty\right\},
\end{equation*}
\begin{equation*}
A\Phi(\gamma):=\sum_{i=1}^\infty \lambda_i(\gamma)  \left\langle S_i,\Phi(\gamma) \right\rangle_{\CM}S_i ,\quad \gamma \in \Omega, \quad \Phi \in D(A).
\end{equation*}
We can then conclude that 
\begin{align}
\label{def_well}
\nonumber \EF (F,F) &= \int \left\langle DF,ADF \right\rangle_\CM \, \varphi d\nu  =\int \left\langle DF,\sum_{i=1}^\infty \lambda_i \left\langle S_i,DF\right\rangle_\CM S_i \right\rangle_{\CM} \, \varphi d\nu \\
& =\sum_{i=1}^\infty \int \lambda_i(\gamma) \left\langle S_i,DF(\gamma)\right\rangle ^2_\CM \, \varphi(\gamma) d\nu < \infty, \quad F\in Y,
\end{align}
is well defined since, for $F\in Y$, this is just a finite sum.

\section{Closability of derivative, divergence, and their inverses}
\label{section_closability}

In this section we are interested in the closability of the derivative, the divergence, and their inverses with respect to a weighted Wiener measure. We use the standard definition of $(\delta,\operatorname{Dom} \delta)$ found in e.g. \cite{Nualart}. For well-definiteness we assume that $\varphi \in L^1(\nu)$. 
The following lemma provides a general density result when using a weighted measure.
\begin{lemma}
\label{lemma_dense_testfunction}
	Let $m$ be a finite Borel measure over a separable metric space $E$ and $\CM$ be a separable Hilbert space.\\
	(a) Let $F\subset L^\infty(E,m)$ be dense in $L^2(E,m/\varphi)$. If $\varphi\in L^1(E,m)$ such that $1/\varphi\in L^1(E,m)$, then $G:=\{f/\varphi:f\in F\}$ is dense in $L^2(E,\varphi m)$.\\
	(b) If $F\subset L^2(E,m)$ is dense in $L^2(E,m)$ then
        $F_\CM:=\left\{\sum_{i=1}^k f\cdot h_i: f\in F, h\in \CM,
        \right.$ $\left.k\in \N\vphantom{\sum_{i=1}^k}\right\}$ is dense in $L^2(E,m;\CM)$.\\
	(c) Let $F$ be as in (a) and $F_\CM$ be defined as in (b). Then $G_\CM:=\left\{\frac{\tilde{f}}{\varphi}: \tilde{f}\in F_\CM\right\}$ is dense in $L^2(E,\varphi m;\CM)$.
	\begin{proof}
		(a) The statement is equivalent to that $\left\{f/\sqrt{\varphi}:f\in F\right\}$ is dense in $L^2(m)$. Approximating any $\Psi\in L^2(m)$ by $f_n/\sqrt{\varphi}$ in $L^2(m)$ is equivalent to approximation of $\Psi\sqrt{\varphi}\in L^2(m/\varphi)$ by $f_n$ in $L^2\left(m/\varphi \right)$. The statement now follows from the density of $F$. \\
		(b) This is a standard result and the proof is therefore omitted.\\
        (c) By (a), $G\subset L^2(E,\varphi m)$ is dense in $L^2(E,\varphi m)$. The statement now follows from (b) with $m$ replaced by $\varphi m$ and $F$ replaced by $G$.
	\end{proof}
\end{lemma}

\begin{remark}
	We note that in Lemma \ref{lemma_dense_testfunction} the condition $F\subset L^\infty(E,m)$ dense in $L^2(E,m/\varphi)$ can be replaced by the condition $F\subset L^\infty(E,m)$ dense in $L^2(E,m)$. 
	\begin{proof}
		We assume $F\subset L^\infty(E,m)$ is dense in $L^2(E,m)$. Let $g\in L^2(E,m/\varphi)$. Given $\eps>0$ let 
		\begin{equation*}
			E^n:=\left\{x:\frac{1}{n}< \frac{1}{\varphi(x)}< n, x\in E\right\}
		\end{equation*}
		and $g_n:=g\cdot \chi_{\{E^n\}}$. It follows $g_n\to g$ in $L^2(m/\varphi)$ when $\nti$ and thus there exists $N\in \N$ such that $\|g_n-g\|_{L^2(m/\varphi)}<\eps/2$ when $n\ge N$. Since $g_n\in L^2(m)$ we can find $f_n\in F$, with support on $E_n$, such that $\|f_n-g_n\|_{L^2(m)}<\eps/2n$ and thus $\|f_n-g_n\|_{L^2(m/\varphi)}<\eps/2$. It follows that $\|f_n-g\|_{L^2(m/\varphi)}<\eps$, $n\ge N$ and the statement is proved.
	\end{proof}
\end{remark}

We now proceed by defining $Z_\CM$ as the set of all cylindrical functions with values in $\CM$ of the form $\sum_{i=1}^k F_j h_j$ where $F_j\in Z$, $h_j\in \CM$ and $k\in \N$. We also let $Z^b$ and $Z_\CM^b$ denote the bounded counterparts to $Z$ and $Z_\CM$ respectively. We define $Z_\CM^p$ as the set of $\CM$-valued cylindrical functions of form $\sum_{i=1}^k \tilde{\phi}_j h_j$ where $h_j\in\CM$, $k\in \N$, $\tilde{\phi}_j$ are polynomials in $\int_0^1H_i\, d\gamma$. Here $H_i,\ i\in \N$, are defined as in \eqref{eq_def_haar1}.
\begin{lemma}
	\label{lemma_zh_dense_1}
	If $x_n\in Z_\CM$, $\CM=L^2(T)$, $T=[0,1]$,
	and
	\begin{equation}
	\label{eq_closability_psidelta}
		\lim_{n \rightarrow \infty} \int \psi \delta(x_n) \, d\nu =0
	\end{equation}
	for all $\psi$ in $Z$. Then 
	\begin{equation*}
		\langle x_n,\tilde{\phi}\rangle_{L^2(\nu;\CM)} \lni 0
	\end{equation*}
	for all $\tilde{\phi}\in Z_\CM^p$.

	\begin{proof}
	For $j\in \N$ and $f\in L^2(T^j)$, symmetric in all $j$ variables, define
	\begin{equation*}
		I_j(f)=j! \int_{t_j=0}^1 \int_{t_{j-1}=0}^{t_j}\ldots \int_{t_1=0}^{t_2}   f(t_1,\ldots,t_j)\, d\gamma_{t_1}\ldots d\gamma_{t_j}.
	\end{equation*}
	We have for all $n\in\N$ the representation
	\begin{equation}
	\label{eq_close_xn}
		x_n(t)=\sum_{j=0}^\infty I_j\big(f_j^{(n)}(\cdot,t)\big),
	\end{equation}
	for some $f_j^{(n)}\in L^2(T^{j+1})$ where for $j\ge 1$, $f_j^{(n)}$ is symmetric in the first $j$ variables and $I_0\big(f_0^{(n)}(t)\big)=E[x_n(t)]$. Thus
	\begin{equation*}
		\delta(x_n)=\sum_{j=0}^\infty I_{j+1}(\tilde{f}_j^{(n)}),
	\end{equation*}
	where $\tilde{f}^{(n)}_j$ denotes the symmetrization of $f_j^{(n)}$. For any $\phi\in \HE_{m+1}$, we have $\phi=I_{m+1}(g)$ for some symmetric $g\in L^2(T^{m+1})$. We even assume that $\phi \in Z$. By hypothesis \eqref{eq_closability_psidelta} we have 
	\begin{align*}
		&(m+1)!\big\langle g,\tilde{f}_m^{(n)}\big\rangle_{L^2(T^{m+1})}
		=\int I_{m+1}(g)I_{m+1}(\tilde{f}_m^{(n)}) \, d\nu\\
		&\qquad=\int \phi\delta(x_n) \, d\nu \rightarrow 0,\text{ as }n\rightarrow \infty.
	\end{align*}
	Since $g\in L^2(T^{m+1})$ is symmetric we even have
	\begin{equation*}
		\big\langle g,f_m^{(n)}\big\rangle_{L^2(T^{m+1})}\rightarrow 0,\text{ as }n\rightarrow \infty.
	\end{equation*}
	For all symmetric polynomials $h\in L^2(T^m;\CM)$ we obtain
	\begin{equation*}
		\langle h,f_m^{(n)}(\cdot,t)\rangle_{L^2(T^m;\CM)}\rightarrow 0,\text{ as }n\rightarrow \infty,
	\end{equation*}
	where the letter $t$ indicates the variable for the function in $\CM$. 
	Using the fact that the $m$th Wiener chaos is isometric to the space of all symmetric functions in $L^2(T^m)$, we get for any $\psi \in L^2(\nu;\CM)$ such that $\psi(\cdot,t)\in \HE_m$ for a.e. $t\in [0,1]$, $\psi=I_m(h)$ that
	\begin{equation*}
		\Big\langle I_m(h),I_m\big(f_m^{(n)}(\cdot,t)\big)\Big\rangle_{L^2(\nu;\CM)}\rightarrow 0,\text{ as }n\rightarrow \infty.
	\end{equation*}
    With \eqref{eq_close_xn} we obtain
    \begin{equation*}
		\langle I_m(h),x_n\rangle_{L^2(\nu;\CM)}\rightarrow 0.
    \end{equation*}
    It follows that
	\begin{equation*}
		\langle x_n,\tilde{\phi}\rangle_{L^2(\nu;\CM)} \lni 0,
	\end{equation*}
for all finite linear combinations
$\tilde{\phi}=\alpha_1\psi_1+\ldots+\alpha_k \psi_k$, and some $k\in
\N$, where the functions $\psi_i\in
L^2(\nu;\CM)$ are in a way that $\psi_i(\cdot,t)$ is from the $i$th Wiener chaos for
a.e. $t\in[0,1]$, and $i\in\{1,\ldots,k\}$. In particular we restrict ourselves to the particular case when $\psi_i$, $i\in\{1,\ldots,k\}$, have a representation of the form 
	\begin{equation}
	\label{eq_close_3}
		\psi_i=\sum_{j=1}^{l}\tilde{\phi}_j(\gamma)\cdot h_j(t),
	\end{equation}
	with $\tilde{\phi}_j\in \HE_m$, $h_j\in \CM$ and some $l\in \N$. Moreover, since polynomials in $L^2([0,1]^m)$ correspond to polynomials in $\HE_m$ we even assume $\tilde{\phi}_i$ to be polynomials in $\int_0^1 H_1\, d\gamma,\ldots,\int_0^1 H_m\,d\gamma$ where $H_1,\ldots,H_m$ are defined in \eqref{eq_def_haar1}. The statement now follows.
	\end{proof}
\end{lemma}

\begin{corollary}
	\label{cor_zh_dense_1}
	Let $x_n\in Z_\CM$ and $x_n\rightarrow x$ in $L^2(\varphi\nu;\CM)$ for some $x\in L^2(\varphi\nu;\CM)$, $\CM=L^2(T)$, $T=[0,1]$. Assume
	\begin{equation}
		\label{closability_cond1}
	\frac{f}{\varphi} \in L^1(\nu),\quad f\in Z,
	\end{equation} 
	and
	\begin{equation}
	\lim_{n \rightarrow \infty} \int \psi \delta(x_n) \, d\nu =0
	\end{equation}
	for all $\psi$ in $Z$. Then $x=0$.
	\begin{proof}
		By \eqref{closability_cond1} we have
		\begin{equation*}
			\lim_{n \rightarrow \infty} \int \left\langle x-x_n,\frac{\tilde{\psi}}{\varphi}\right\rangle_\CM  \, \varphi d\nu =0,
		\end{equation*}
		for all $\tilde{\psi}\in Z_\CM$, thus 
		\begin{equation*}
			\lim_{n \rightarrow \infty} \int \left\langle x_n,\tilde{\psi}\right\rangle_\CM \,  d\nu =\int \left\langle x,\tilde{\psi}\right\rangle_\CM \,  d\nu,\quad \tilde{\psi}\in Z_\CM. 
		\end{equation*}
		By Lemma \ref{lemma_zh_dense_1} the right hand side is $0$ for all $\tilde{\psi}\in Z_\CM^p$. Since $Z_\CM^p$ is dense in $L^2(\nu;\CM)$ the statement follows.
	\end{proof}
\end{corollary}

We emphasize that Corollary \ref{cor_zh_dense_1}, which is based on Lemma \ref{lemma_zh_dense_1}, is crucial for the proof of Proposition \ref{prop:close_delta_inv} below. In addition we get the following proposition as an independent result.

\begin{proposition}
If \eqref{closability_cond1} holds then $\{\frac{1}{\varphi}D\psi:\psi \in Z\}$ is a dense subset of $L^2(\varphi \nu;\CM)$.
\begin{proof}
We assume the contrary. Then we can find $x\in L^2(\varphi \nu;\CM)$ and $x_n\in Z_\CM$, $n\in\N$, with values in $\CM$, such that $x\neq 0$ and $x_n\lni x$ in $L^2(\varphi \nu;\CM)$ and
\begin{equation*}
\int \frac{1}{\varphi}\langle D\psi,x\rangle_\CM \varphi \, d\nu=0
\end{equation*}
for all $\psi \in Z$. It follows that
\begin{align*}
0&=\int\left\langle \frac{D\psi}{\varphi},x\right\rangle_\CM  \,\varphi d\nu
=\lim_{n\to \infty} \int \left\langle \frac{D\psi}{\varphi},x_n\right\rangle_\CM
\, \varphi d\nu
=\lim_{n\to \infty} \int \left\langle D\psi,x_n\right\rangle_\CM
\,  d\nu\\
&=\lim_{n\to \infty} \int \psi \delta(x_n)  \, d\nu
\end{align*}
for all $\psi \in Z$. From this we get $x=0$ by Corollary \ref{cor_zh_dense_1}, and we have a contradiction.
\end{proof}
\end{proposition}

\begin{proposition}
\label{lemma_closable1}
	Assume $\varphi\in L^1(\nu)$ and \eqref{closability_cond1} holds, i.e.
	\begin{equation*}
		\frac{f}{\varphi} \in L^1(\nu),\quad f\in Z.
	\end{equation*}
	Then $(\mathscr{D},Z)$ defined by
	\begin{equation*}
		\mathscr{D}(f,g):=\hf \int \langle Df,Dg\rangle_\CM \, \varphi d\nu, \quad f,g\in Z,
	\end{equation*}
	is closable on $L^2(\varphi \nu)$ and we denote this closure by $\D^1$. 
	We also note that $(D,Z)$ is closable as an operator $L^2(\varphi\nu)\supset Z \to L^2(\varphi\nu;\CM)$.

	\begin{proof}
	We observe that under \eqref{closability_cond1}, $\delta(\psi)/\varphi\in L^2(\varphi\nu)$ for $\psi\in Z_\CM^b:=\{\sum_{i=1}^k \phi_ih_i: k\in \N, \phi_i\in Z\cap L^\infty(\nu), h_i\in\CM\}$. Now for $u_n\in Z$, $\psi \in Z_\CM$ and
	\begin{equation*}
		u_n \lni 0 \quad \text{in } L^2(\varphi \nu), \quad Du_n 
		\lni f \quad \text{in } L^2(\varphi \nu; \CM),
	\end{equation*}
	we have
	\begin{align}
	\label{eq_closability_D_un}
		\hf \int \left\langle Du_n,\frac{\psi}{\varphi} \right\rangle_\CM \, \varphi d\nu			
		&=\hf \int  u_n\frac{\delta(\psi)}{\varphi}  \, \varphi d\nu.
	\end{align}
By the earlier observation the right hand side of \eqref{eq_closability_D_un} tends to $0$. By Lemma Lemma \ref{lemma_dense_testfunction}(c), test functions of form $\psi/\varphi$, $\psi\in Z^b_\CM$ are dense in $L^2(\varphi\nu;\CM)$. It follows that $f=0$ and thus $\mathscr{D}$ is closable on $L^2(\varphi\nu)$. 
 	\end{proof}
\end{proposition}

\begin{proposition} 
\label{prop:close_delta}
 	If \eqref{closability_cond1} holds then	$(\delta,Z_\CM)$ is a closable operator in $L^2(\varphi\nu;\CM)\supset Z_\CM \to L^2(\varphi\nu)$.
 	\begin{proof}
 		Let $x_n\in Z_\CM$ such that 
		\begin{align*}
		&x_n\lni 0 \quad \text{in }L^2(\varphi\nu),\quad \delta(x_n)\lni f\quad \text{in }L^2(\varphi\nu;\CM).
		\end{align*}
		Then for all $y\in \{z/\varphi: z\in Z^b\}$ we have
		\begin{align*}
			\langle f,y\rangle_{L^2(\varphi\nu)}
			&=\lim_\nti \langle \delta(x_n),y\rangle_{L^2(\varphi\nu)}\\
			&=\lim_\nti \left\langle x_n,\frac{D(y\varphi)}{\varphi}\right\rangle_{L^2(\varphi\nu;\CM)}=0.
		\end{align*}
 	\end{proof}
 \end{proposition}

We proceed by defining equivalence classes in $L^2(\varphi\nu)$. Let $x\in L^2(\varphi\nu)$ be the representative of $\Bx(x):=\{x+c\I: c\in  \R\},\ x\in L^2(\varphi\nu)$ where $\I$ is the constant function taking the value $1$ on $\Omega$. We also define
\begin{equation*}
	L^2(\varphi\nu)\ominus \I:=\{\Bx(x):x\in L^2(\varphi\nu)\}
\end{equation*}
as well as $\Bdelta(z)$ which denotes the equivalence class of $\delta(z)$, $z\in Z_\CM$. We remark that every $\Bx\in L^2(\varphi\nu)\ominus \I$ has a unique representative $x_0\equiv x_0(\Bx)$ with $\int x_0 \, \varphi d\nu=0$. $L^2(\varphi\nu)\ominus \I$ becomes a Hilbert space with the inner product
\begin{equation*}
\langle \Bx_1,\Bx_2\rangle_{L^2(\varphi\nu)\ominus \I}:=\langle x_0(\Bx_1),x_0(\Bx_2)\rangle_{L^2(\varphi\nu)},\quad \Bx_1,\Bx_2\in L^2(\varphi\nu)\ominus \I.
\end{equation*}

\begin{proposition}
\label{prop:close_delta_inv}
	Assume that \eqref{closability_cond1} holds. Then $\delta^{-1}$ defined on $Z_\CM^{-1}:=\{\Bdelta(z):z\in \Z_\CM\}$ as an operator $L^2(\varphi\nu)\ominus \I \supset Z_\CM^{-1} \to L^2(\varphi\nu;\CM)$ is closable.
	\begin{proof}
		Let $y_n\in Z_\CM^{-1}$ such that 
		\begin{align*}
		&y_n\lni 0 \quad \text{in }L^2(\varphi\nu)\ominus \I,\quad \delta^{-1}(y_n)\lni f\quad \text{in }L^2(\varphi\nu;\CM).
		\end{align*}
		Let $y_n=\Bdelta(z_n)$ for some $z_n\in Z_\CM$. Then 
		\begin{align*}
			\Bdelta(z_n)\lni 0 \quad \text{in }L^2(\varphi\nu)\ominus \I
		\end{align*}
		which is equivalent to $x_0(\Bdelta(z_n))\lni 0$ in $L^2(\varphi\nu)$. Noting that $x_0(\Bdelta(z_n))=\delta(z_n)+c_n\I$ for some $c_n\in \R$. We have
		\begin{align*}
			c_n
			&=\langle \delta(z_n)+c_n\I,\I\rangle_{L^2(\nu)}
			=\left\langle \delta(z_n)+c_n\I,\frac{\I}{\varphi}\right\rangle_{L^2(\varphi\nu)}\\
			&=\left\langle x_0(\Bdelta(z_n)),\frac{\I}{\varphi}\right\rangle_{L^2(\varphi\nu)}\lni 0.
		\end{align*}
		Summarizing we get
		\begin{align*}
			 \delta(z_n)\lni 0 \quad \text{in }L^2(\varphi\nu)\quad \mbox{\rm and}\quad  z_n=\delta^{-1}(\Bdelta(z_n))\lni f \quad \quad\text{in }L^2(\varphi\nu;\CM),
		\end{align*}
		and conclude
		\begin{align*}
			\int \psi \delta(z_n)=\left\langle \frac{\psi}{\varphi},\delta(z_n)\right\rangle_{L^2(\varphi\nu)} \lni 0
		\end{align*}
		for all $\psi \in Z$. By Corollary \ref{cor_zh_dense_1} we have $f=0$. The statement follows.
	\end{proof}
\end{proposition}

In order to have a well defined inverse of the gradient operator, we define equivalence classes in $Z$,
\begin{equation*}
	Z^e:=\{\Bx(x):x\in Z\},
\end{equation*}
and
\begin{equation*}
	D^e\Bx\equiv D^e\Bx(x):=Dx, \quad \Bx\in Z^e.
\end{equation*}
$D^e$ is invertible on $Z^e$. We denote the inverse by $D^{e-1}$ and the corresponding domain by $Z^{e-1}:=\{D^e(\Bx):\Bx\in Z^e\}$.

\begin{proposition}
\label{prop:close_gradient}
	Assume \eqref{closability_cond1}. The operator $(D^{e-1},Z^{e-1})$ as an operator $L^2(\varphi\nu;\CM)\supset Z^{e-1}\to L^2(\varphi\nu)\ominus \I$ is closable.
	\begin{proof}
		Let $x_n\in Z^{e-1}$ such that
		\begin{align*}
			&x_n\lni 0 \quad \text{in }L^2(\varphi\nu;\CM),\quad D^{e-1}(x_n)\lni f\quad \text{in }L^2(\varphi\nu)\ominus \I.
		\end{align*}
		For $x_n=D^e\Bx_n$, $\Bx_n\in Z^e$, and all $\By\in \tilde{Y}:=\{\Bx(y):y\varphi=\delta(z),\text{ for some }z\in Z_\CM^b\}$, we have
		\begin{align}
		\notag
			\langle f,\By\rangle_{L^2(\varphi\nu)\ominus \I}
			\notag
			&=\lim_\nti \left\langle x_0(D^{e-1}x_n),x_0(\By)\right\rangle_{L^2(\varphi\nu)}\\
			\notag
			&=\lim_\nti \left\langle x_0(\Bx_n),x_0(\By)\varphi\right\rangle_{L^2(\nu)}\\
			\notag
			&=\lim_\nti \left\langle x_0(\Bx_n),\delta(z)\right\rangle_{L^2(\nu)}\\
			\notag
			&=\lim_\nti \left\langle D^e(\Bx_n),z\right\rangle_{L^2(\nu;\CM)}\\
			\label{eq:dinv:1}
			&=\lim_\nti \left\langle x_n,\frac{z}{\varphi}\right\rangle_{L^2(\varphi\nu;\CM)}=0.
		\end{align}
		We now show that $\tilde{Y}$ is dense in $L^2(\varphi\nu)\ominus\I$. Since $0=E\left[\frac{f}{\varphi}{\delta(z)}\right]$, $z\in Z_\CM^b$ implies $f=c\cdot\varphi$, $c\in \R$, we obtain that $\left\{\delta(z)+c\I:c\in \R, z\in Z_\CM^b\right\}$ is dense in $L^2(\nu)$. Now Lemma \ref{lemma_dense_testfunction} says that $\left\{\frac{\delta(z)+c\I}{\varphi}: c\in \R,z\in Z_\CM^b\right\}$ is dense in $L^2(\varphi\nu)$. Using the fact that
		\begin{equation*}
			\int \frac{\delta(z)+c\I}{\varphi}\, \varphi d\nu=0, \quad z\in Z_\CM^b
		\end{equation*}
		yields that $c=0$, it follows that $\left\{\frac{\delta(z)}{\varphi}: z\in Z_\CM^b\right\}$ is dense in $L^2(\varphi\nu)\ominus \{c\I:c\in \R\}$. Thus $\tilde{Y}$ is dense in $L^2(\varphi\nu)\ominus\I$. \eqref{eq:dinv:1} shows that $f=0$ and the claim follows.
	\end{proof}
\end{proposition}

\begin{remark}
Under the condition \eqref{closability_cond1}, i.e. 
\begin{equation*}
	\frac{f}{\varphi} \in L^1(\nu),\quad f\in Z,
\end{equation*}
the operators $D$, $\delta$, $D^{e-1}$ and $\delta^{-1}$ are closable in the sense of Propositions \ref{lemma_closable1}, \ref{prop:close_delta}, \ref{prop:close_delta_inv}, and \ref{prop:close_gradient}.
\end{remark}
In the following sections we are interested in weight functions of form 
	\begin{align}
		\label{eq:exp:phi}
		\varphi(\gamma):=\exp\left\{\int_0^1\langle b(\gamma_s),d\gamma_s\rangle_{\R^d}-\hf\int_0^1|b(\gamma_s)|^2\, ds\right\},
	\end{align}
	where $\int_0^1\langle b(\gamma_s),d\gamma_s\rangle_{\R^d}$ is a stochastic It\^o integral in the sense of \cite{RevuzYor}. It is natural to ask for conditions on $b$ so that $\frac{f}{\varphi}\in L^1(\nu)$, $f\in Z$. Related questions have been addressed in great detail in \cite{Stummer1993} which we follow for the subsequent results.
	
\begin{proposition}
	Let $d\ge 3$ and $\varphi$ be given by \eqref{eq:exp:phi} where $b=(b^1,\ldots,b^d)$ such that $b^j$ is for all $j=1,\ldots,d$ continuous of type $b^j:\R^d\to \overline{\R}$, $\overline{\R}$ being endowed with a suitable topology.	\\
	(a) Let $b$ be spherically symmetric, i.e. $b(x)=g(|x|)$, $x\in \R^d$, for some $g:\R\to \R^d$. Then $\frac{f}{\varphi}\in L^1(\nu)$ for all $f\in Z$ if and only if  
	\begin{equation}
		\label{eq:exp:gr}
		\int_0^1 rg(r)^2\, dr<\infty\quad\text{and}\quad\sup_{s\ge 2} \int_{s-1}^{s+1} g(r)^2 \, dr<\infty.
	\end{equation}
	(b) In the case of a non-spherically symmetric $b$,
	\begin{align}
		\label{eq:exp:5}
		\sup_{x\in \R^d} \int_{\|x-y\|< 1} |b(y)|^p \, dy<\infty\quad\text{for some }p>d,
	\end{align}
	implies  $\frac{1}{\varphi^q}\in L^1(\nu)$ for all $q\ge 1$ and hence $\frac{f}{\varphi}\in L^1(\nu)$ for all $f\in Z$. 
\end{proposition}

\begin{proof}
	In order to follow \cite{Stummer1993} it is important to note that the process $X$ in this reference is a Brownian motion with respect to the measure $P_\cdot$.\\
	(a) \textit{Step 1:} We show sufficiency of \eqref{eq:exp:gr}. Let 
	\begin{align*}
		Y(t):=\exp\left\{-2\int_0^t\langle b(\gamma_s),d\gamma_s\rangle_{\R^d}-\int_0^t|b(\gamma_s)|^2\, ds\right\},\quad t\in [0,1].
	\end{align*}
	Using the Cauchy--Schwartz inequality we obtain
	\begin{align*}
		E\left[\frac{f}{\varphi}\right]
		&=E\left[f(\gamma)\cdot \exp\left\{\hf\int_0^1|b(\gamma_s)|^2\, ds-\int_0^1\langle b(\gamma_s),d\gamma_s\rangle_{\R^d}\right\}\right]^2\\
		&\le E\left[f(\gamma)^2\cdot \exp\left\{2\int_0^1|b(\gamma_s)|^2\, ds\right\} \right]\cdot E\left[ Y(1)\right].
	\end{align*}
	Here the second factor is bounded by $1$ since $Y(t)$, $t\in[0,1]$, is a supermartingale, see \cite{KaratzasShreve} chapter 3.5 D. Noting that $f^4\in L^1(\nu)$, it remains to show that
	\begin{align*}
		E\left[\exp\left\{4\int_0^1|b(\gamma_s)|^2\, ds\right\} \right]<\infty.
	\end{align*}
	However, this is an immediate consequence of Theorem 2.4 in \cite{Stummer1993}. \\
	\textit{Step 2:} We show that \eqref{eq:exp:gr} is necessary. Let
	\begin{align*}
		I(\gamma):=\int_0^1\langle b(\gamma_s),d\gamma_s\rangle_{\R^d}.
	\end{align*}
	According to the above hypotheses on $b$, we may assume well-definiteness of $I(\cdot)$ as a random variable. In addition, because of the spherical symmetry of $b$ we may suppose that 
	$I(\gamma)=-I(-\gamma)$. Since $1/\varphi \in L^1(\nu)$ we have
	\begin{align*}
		\infty 
		&>E\left[\exp\left\{\hf\int_0^1|b(\gamma_s)|^2\, ds-\int_0^1\langle b(\gamma_s),d\gamma_s\rangle_{\R^d}\right\}\right] \\
		&\ge E\left[\chi_{\left\{I(\gamma)< 0\right\}}\cdot\exp\left\{\hf\int_0^1|b(\gamma_s)|^2\, ds-\int_0^1\langle b(\gamma_s),d\gamma_s\rangle_{\R^d}\right\}\right] \\
		&\ge E\left[\chi_{\left\{I(\gamma)< 0\right\}}\cdot \left(\hf\int_0^1|b(\gamma_s)|^2\, ds-\int_0^1\langle b(\gamma_s),d\gamma_s\rangle_{\R^d}\right)\right] \\
		&\ge E\left[\chi_{\left\{I(\gamma)> 0\right\}}\cdot \int_0^1\langle b(\gamma_s),d\gamma_s\rangle_{\R^d}\right],
	\end{align*}
	with $\chi$ denoting the indicator function, and thus 
	\begin{align}
		\label{eq:exp:bounded}
		\left|E\left[\int_0^1\langle b(\gamma_s),d\gamma_s\rangle_{\R^d}\right]\right|<\infty.
	\end{align}
	
	Using Jensen's inequality and recalling \eqref{eq:exp:phi}, it follows from $1/\varphi \in L^1(\nu)$ that
	\begin{align}
		\label{eq:exp:jensen1}
		\exp\left\{E\left[\hf\int_0^1|b(\gamma_s)|^2\, ds-\int_0^1\langle b(\gamma_s),d\gamma_s\rangle_{\R^d}\right]\right\}<\infty.
	\end{align}
	From \eqref{eq:exp:jensen1} and \eqref{eq:exp:bounded} we obtain
	\begin{align}
		\label{eq:exp:1}
		E\left[\int_0^1|b(\gamma_s)|^2\, ds\right]<\infty.
	\end{align}
	Since $\gamma$ is a $d$-dimensional Wiener process and $|b(\gamma_s)|^2=b(|\gamma_s|)^2$ we get
	\begin{align}
		\label{eq:exp:2}
		E\left[\int_0^1|b(\gamma_s)|^2\, ds\right]
		= \int_{r=0}^\infty b(r)^2\int_{s=0}^1 \int_{y\in S_r} (2\pi s)^{-d/2}\cdot e^{-r^2/2s} \,dy\,ds\,dr,
	\end{align}
	where $S_r$ denotes the $(d-1)$-sphere of radius $r$. Equations \eqref{eq:exp:1} and \eqref{eq:exp:2} give
	\begin{align*}
		\infty
		&>\int_{r=0}^1 b(r)^2\int_{s=0}^1 \frac{r^{d-1}}{s^{d/2}}e^{-r^2/2s} \,ds\,dr+\int_{r=1}^\infty b(r)^2\int_{s=0}^1 \frac{r^{d-1}}{s^{d/2}}e^{-r^2/2s} \,ds\,dr \\
		&=C\int_{r=0}^1 b(r)^2r\,dr+C\int_{r=1}^\infty b(r)^2r\,dr\\
		&\ge C\int_{r=0}^1 b(r)^2r\,dr+C\sup_{s\ge 2}\int_{s-1}^{s+1} b(r)^2\,dr,
	\end{align*}
	where the constant $C$ results from integration with respect to $s$. The statement follows.
	\\(b) We follows exclusively the argumentation of \cite{Stummer1993}. Our relation \eqref{eq:exp:5} is (1.12) in \cite{Stummer1993}. This implies condition (1.5) of \cite{Stummer1993} by the proof of Corollary 1.7 in \cite{Stummer1993}. Using Theorem 1.2 of \cite{Stummer1993} we get (1.1) of \cite{Stummer1993}. Replacing $b$ by $-b$ this yields $\frac{1}{\varphi^q}\in L^1(\nu)$ for all $q\ge 1$  according to Theorem 5.1 of \cite{Stummer1993}.
\end{proof}

\section{Closability of the bilinear form}
\label{section_bilinear}

We formulate the following lemma that serves an important role in proving and formulating the closability conditions of the form.
\begin{lemma}
\label{lemma1}
If $s_k=\sum_{i=1}^{r}c_i\cdot2^{-i}$, $j\in \{1,\dots,d\}$ and
\begin{equation*}
	i_{p,j}\equiv i_{p,j}(c_0,c_1,\ldots):=
	\begin{cases}
		2^{p-1}(d+j-1)+1+ \sum_{q=0}^{p-1} c_q 2^{p-q-1}, &p\ge 1, \\
		j, &p=0, 
	\end{cases}
\end{equation*}
where $c_1,\dots,c_{r-1} \in \{0,1\}$, $c_r=1$, $c_0=c_{r+1}=0$ and $r\ge 2$, then
\begin{equation*}
\label{eq_lambda_sum}
\sum_{i=1}^{\infty} \lambda_i \langle S_i(s_k), e_j\rangle_{\RE^d}^2=\lambda_j s_k^2+\sum_{p=1}^{r}\lambda_{i_{p,j}}2^{p-1}\left(c_{p}2^{-p}+(-1)^{c_p}\sum_{q=p+1}^{r+1 } 2^{-q}c_q\right)^2.
\end{equation*}
(b) The relation
\begin{equation}
\label{eq_lambda_old_bound}
\sup_{\begin{smallmatrix} c_1,\dots,c_{r-1} \in \{0,1\}, \\ c_r=1,c_{r+1}=0, \\r\ge 2 \end{smallmatrix}}\int \sum_{j=1}^d \sum_{p=1}^{r}\lambda_{i_{p,j}}2^{p-1}\left(c_{p}2^{-p}+(-1)^{c_{p}}\sum_{q=p+1}^{r+1} 2^{-q}c_q\right)^2\, \varphi d\nu< \infty
\end{equation}
is equivalent to 
\begin{equation}
\label{eq_lambda_new_bound}
\sup_{c_1,c_2,\ldots \in \{0,1\}}\sum_{j=1}^d \sum_{p=1}^\infty \frac{\int\lambda_{i_{p,j}}\,\varphi d\nu}{2^p} < \infty.
\end{equation}

\begin{proof}
Since $e_j$ $\bot$ $S_i$ for $i\mod d \neq j$ and $S_i(s_k)=0$ for $i>d2^r$ we get 
\begin{align*}
\sum_{i=1}^{\infty} \lambda_i \langle S_i(s_k), e_j\rangle_{\R^d}^2&=\lambda_j \langle S_j(s_k), e_j\rangle_{\R^d}^2+\sum_{i=1}^{2^r} \lambda_{di+j} \langle S_{di+j}(s_k), e_j\rangle_{\R^d}^2.
\end{align*}
Now since
\begin{equation*}
\lambda_j\langle S_j(s_k), e_j\rangle_{\R^d}^2=\lambda_j s_k^2,
\end{equation*}
and 
\begin{equation*}
S_{di+j}(s_k)=0
\end{equation*} 
unless $di+j=2^{p-1}(d+j-1)+1+ \sum_{q=0}^{p-1} c_q 2^{p-q-1}$ for some $p\ge 1$. We get
\begin{align*}
&\sum_{i=1}^{2^r} \lambda_{di+j} \langle S_{di+j}(s_k), e_j\rangle_{\R^d}^2\\
&\quad=\sum_{p=1}^{r} \lambda_{i_{p,j}} \langle S_{i_{p,j}}(s_k), e_j\rangle_{\R^d}^2,
\end{align*}
where
\begin{align*}
\langle S_{i_{p,j}}(s_k), e_j\rangle_{\R^d}&=\begin{cases}
2^{\frac{p-1}{2}}\sum_{q=p+1}^{r+1} 2^{-q}c_q &\mbox{if } \text{$c_{p}=0$} \\
2^{-\frac{p+1}{2}}-2^{\frac{p-1}{2}}\sum_{q=p+1}^{r+1} 2^{-q}c_q &\mbox{if } \text{$c_{p}=1$}
\end{cases} \\ &=2^{\frac{p-1}{2}}\Big(c_p 2^{-p}+(-1)^{c_p}\sum_{q=p+1}^{r+1} 2^{-q}c_q\Big),
\end{align*}
and thus the statement follows.
\\(b) Assuming \eqref{eq_lambda_old_bound} and recalling that $\lambda>0$. 
\begin{align*}
\begin{split}
	&\sup_{\begin{smallmatrix} c_1,\dots,c_{r-1} \in \{0,1\}, \\ c_r=1,c_{r+1}=0,\\ j\in \{1,\dots,d\}, r\ge 2 \end{smallmatrix}} \int\sum_{j=1}^d\sum_{p=1}^{r}\lambda_{i_{p,j}}2^{p-1}\left(c_{p}2^{-p}+(-1)^{c_{p}}\sum_{q=p+1}^{r+1} 2^{-q}c_q\right)^2 \,\varphi d\nu\\
	&\qquad \ge \sup_{c_1,c_2,\ldots \in \{0,1\}}\int \sum_{j=1}^d\sum_{p=1}^{\infty}\lambda_{i_{p,j}}2^{p-1}\left(2^{-p-1} \right)^2\, \varphi d\nu \\
	&\qquad= \sup_{c_1,c_2,\ldots \in \{0,1\}} \frac{1}{8} \sum_{j=1}^d\sum_{p=1}^{\infty}\frac{\int\lambda_{i_{p,j}}\,\varphi d\nu}{2^{p}}.
\end{split}
\end{align*}
Assuming \eqref{eq_lambda_new_bound} we obtain
\begin{align*}
\begin{split}
	& \sup_{\begin{smallmatrix} c_1,\dots,c_{r-1} \in \{0,1\}, \\ c_r=1,c_{r+1}=0,\\ r\ge 2 \end{smallmatrix}} \int \sum_{j=1}^d \sum_{p=1}^{r}\lambda_{i_{p,j}}2^{p-1}\left(c_{p}2^{-p}+(-1)^{c_{p}}\sum_{q=p+1}^{r+1} 2^{-q}c_q\right)^2 \,\varphi d\nu\\
	&\qquad \le \sup_{c_1,c_2,\ldots \in \{0,1\}} \int \sum_{j=1}^d\sum_{p=1}^{\infty}\lambda_{i_{p,j}}2^{p-1}\left(2^{-p} \right)^2 \,\varphi d\nu \\
	&\qquad = \sup_{c_1,c_2,\ldots \in \{0,1\}} \hf \sum_{j=1}^d\sum_{p=1}^{\infty}\frac{\int\lambda_{i_{p,j}}\, \varphi d\nu}{2^{p}}.
\end{split}
\end{align*}
The statement now follows.
\end{proof}
\end{lemma}

\begin{remark}
	We are most interested in investigating the degree of increase for the eigenvalues. Keeping this in mind, it makes sense to assume that the sequence of eigenvalues is non-decreasing $0<\lambda_1(\gamma)<\lambda_2(\gamma),\ldots$. Under this condition, relation \eqref{eq_lambda_new_bound} simplifies to
	\begin{equation*}
		\sum_{p=0}^\infty  \frac{\int \lambda_{d2^p} \, \varphi d\nu}{2^p} < \infty.
	\end{equation*}
\end{remark}

\begin{proposition}
\label{prop_close_1}
Let $\varphi$ satisfy \eqref{closability_cond1} and assume that $A\ge \varepsilon \operatorname{Id}$ for some $\varepsilon>0$.\\
(a) The form $(\EF,Y)$ is closable in $L^2(\varphi \nu)$. Let $(\EF,D_Y(\EF))$ denote the closure of $(\EF,Y)$ in $L^2(\varphi \nu)$.
\\(b) If 
\begin{equation}
\label{quasi_cond}
\sup_{c_1,c_2,\ldots \in \{0,1\}} \sum_{j=1}^d \sum_{p=0}^\infty \frac{\int\lambda_{i_{p,j}}\,\varphi d\nu}{2^p} < \infty,
\end{equation}
where $i_{p,v}$ and $c_1,c_2,\ldots$ are constructed as in Lemma \ref{lemma1}, then $(\EF,Z)$ is closable and we let $(\EF,D_Z(\EF))$ denote the closure of $(\EF,Z)$ on $L^2(\varphi \nu)$. Furthermore $D_Z(\EF)=D_Y(\EF)$ under \eqref{quasi_cond}.
\\(c) We have $Z\subset D_Y(\EF)$ if and only if (\ref{quasi_cond}) holds. 
\\(d) If $(\EF,Z)$ is closable in $L^2(\varphi \nu)$ then \eqref{quasi_cond} holds. 
\begin{proof}
(a) We have
\begin{align*}
\mathcal{E}(F,F)&=\sum_{i=1}^{\infty}  \int \lambda_i \left\langle S_i,DF\right\rangle_\mathbb{H}^2  \, \varphi d\nu=\int\left\langle A^{1/2}DF,A^{1/2}DF \right\rangle_\mathbb{H} \, \varphi d\nu.
\end{align*}
Suppose $\{F_n\}_{n\ge 1}\subset Y$ such that $F_n \lni 0$ in $L^2(\varphi \nu)$ and $\mathcal{E}(F_n-F_m,F_n-F_m)\rightarrow 0$. Since it follows that $A^{1/2}DF_n$ is Cauchy in $L^2(\varphi \nu;\CM)$ we may define
\begin{equation*}
\psi:=\lim_{n\rightarrow \infty} A^{1/2}DF_n.
\end{equation*}
We also define 
\begin{equation}
\label{eq3}
JG:=\sum_{i=1}^{\infty}\lambda_i^{-1/2} \left\langle S_i,G\right\rangle_\mathbb{H} S_i\quad G\in L^2(\varphi\nu;\CM).
\end{equation}
The operator $J$ is bounded on $L^2(\varphi \nu;\CM)$ and it follows
\begin{equation*}
DF_n=JA^{1/2}DF_n \lni J\psi \text{ in } L^2(\varphi \nu;\CM).
\end{equation*}
	From Proposition \ref{lemma_closable1}, it is known that $(\mathscr{D},Z)$ is closable on $L^2(\varphi \nu)$. It follows that $DF_n\lni 0$ and thus $J\psi=0$. Since $\lambda_i>0$ and $\lambda_i^{-1/2}>0$, (\ref{eq3}) gives $\psi=0$. Now $A^{1/2}DF_n\lni  0$ and thus $\mathcal{E}(F_n,F_n)=\int\left\langle A^{1/2}DF_n,A^{1/2}DF_n \right\rangle_\mathbb{H} \, \varphi d\nu \rightarrow  0$ as $n \rightarrow \infty$.
	\\(b) We show that (\ref{quasi_cond}) implies $Z\subset D_Y(\EF)$. Then $Y\subset Z\subset D_Y(\EF)$. Since $(\EF,Y)$ is closable with closure $(\EF,D_Y(\EF))$, cf. (a), $(\EF,Z)$ is then also closable and has $(\EF,D_Y(\EF))$ as its closure, i.e. $D_Y(\EF)=D_Z(\EF)$. 
	\\
	Let $x^v(p)$ denote the $v$th coordinate of $p\in  \RE^d$, $v\in \{1,\dots,d\}$. Let us demonstrate that $F(\gamma)=x^v(\gamma(s))\in D_Y(\EF)$ for all $v\in \{1,\dots,d\}$ and all $s\in [0,1]$. Fix $s\in [0,1]$, $v\in \{1,\dots,d\}$, and let $s_k \lki s$ where $s_k$ is a sequence of dyadic numbers. Now let
	\begin{equation*}
		F_{v,k}(\gamma):=x^v(\gamma(s_k))\in Y \subset D_Y(\EF),\quad k\in \N.
	\end{equation*}
	We aim to apply Lemma \ref{lemma1}. By \eqref{quasi_cond} we have 
	\begin{align}
	\nonumber 
		&\int\sum_{i=1}^\infty \lambda_i(\gamma)\langle S_i,DF(\gamma)\rangle_\CM^2 \,\varphi d\nu
		=\int \sum_{i=1}^\infty \lambda_i(\gamma)\langle
                S_i(s),e_v\rangle_{\RE^d}^2 \,\varphi d\nu \\
        \nonumber
		&\qquad \le \sup_k \int \sum_{i=1}^\infty \lambda_i(\gamma)\langle S_i(s_k),e_v\rangle_{\RE^d}^2 \, \varphi d\nu \\
		\label{ubound}
		&\qquad \le \sum_{v=1}^d\int \sum_{p=0}^\infty \frac{\lambda_{i_{p,v}}}{2^p} \, \varphi d\nu
		<\infty.
	\end{align}
	Furthermore $F_{v,k} \lki F$ in $L^2(\varphi \nu)$ and
	\begin{equation*}
		\int \sum_{i=1}^\infty \lambda_i(\gamma)\langle S_i,D(F-F_{v,k})(\gamma)\rangle_\CM^2 \, \varphi d\nu
		\le \int \sum_{i=n_k}^\infty \lambda_i(\gamma)\langle S_i,DF(\gamma)\rangle_\CM^2 \, \varphi d\nu
	\end{equation*}
	for some sequence $n_k\rightarrow \infty$ as $k\rightarrow \infty$. This says $F_{v,k} \lki F$ in $\EF_1$-norm and $F\in D_Y(\EF)$, where we recall that $\|\cdot\|_{\EF_1}^2=\EF(\cdot,\cdot)+\|\cdot\|^2_{L^2(\varphi \nu)}$. Indeed, the same conclusion applies to an arbitrary $F \in Z$ with $F(\gamma)=f\left(\gamma(s_1),\dots,\gamma(s_k)\right)$, $s_1,\dots,s_k\in[0,1]$ since because of \eqref{def_grad} we have
	\begin{equation*}
	\langle S_i,DF(\gamma)\rangle_\CM=\sum_{j=1}^d \sum_{i'=1}^k \frac{\partial f}{\partial x_{i',j}}\left(\gamma(s_1),\dots,\gamma(s_k)\right)\langle S_i(s_{i'}),e_j\rangle_{\RE^d}.
	\end{equation*}
	We then apply \eqref{def_well} and \eqref{ubound}.
	\\(c) Recalling (b), we still have to show that $Z \subset D_Y(\EF)$ implies (\ref{quasi_cond}). 
	\\
	We suppose $Z\subset D_Y(\EF)$. Therefore $x_v(\gamma(s))\in D_Y(\EF)$ for all $v\in \{1,\dots,d\}$ and all $s\in [0,1]$. First we check the expression of Lemma \ref{lemma1} for all dyadic points $s_k$. To get a bound on 
	\begin{equation*}
		\int\sum_{i=1}^\infty \lambda_i\langle S_i(s), e_v\rangle_{\RE^d}^2\,\varphi d\nu,
	\end{equation*}
	when $s\in [0,1]$ is no longer dyadic, we need \eqref{quasi_cond}. The statement follows.
		\\(d) $(Z,\EF)$ closable in $L^2(\varphi\nu)$ implies $\EF_1(F,F)<\infty$ for all $F\in Z$, in particular for $F(\gamma)=x^v(\gamma(s))$ for all $v\in\{1,\ldots,d\}$ and all $s\in[0,1]$. We have
		\begin{equation*}
			\EF_1(F,F)=\|F\|^2_{L^2(\varphi\nu)}+\int\sum_{i=1}^\infty \lambda_i \langle S_i(s),e_v\rangle^2_{\RE^d}\, \varphi d\nu<\infty.
		\end{equation*}
		Thus $\int \sum_{i=1}^\infty \lambda_i \langle S_i(s),e_v\rangle^2_{\RE^d}\, \varphi d\nu<\infty$ for all $v\in\{1,\ldots,d\}$ and all $s\in[0,1]$. Now Lemma \ref{lemma1}(a) and (b) implies \eqref{quasi_cond}.
\end{proof}
\end{proposition}

\begin{proposition}
Let $\varphi$ satisfy \eqref{closability_cond1} and \eqref{quasi_cond}, and assume $A\ge \varepsilon \operatorname{Id}$ for some $\varepsilon>0$. Then the form $(\EF,D_Z(\EF))$ is a Dirichlet form on $L^2(\varphi \nu)$.

	\begin{proof}
	We use Proposition I$.4.10$ from \cite{Rockner}. It follows that we must show $\mathcal{E}(1\wedge F^+,1\wedge F^+)\le \mathcal{E}(F,F)$. We know that for $F\in Y$
	\begin{align*}
		\mathcal{E}(F,F) &=\sum_{i=1}^{\infty} \int \lambda_i\left\langle S_i,DF\right\rangle_\mathbb{H}^2  \, \varphi d\nu
		=\sum_{i=1}^{\infty} \int  \lambda_i (\dSi,F)^2  \, \varphi d\nu \\
		&=\sum_{i=1}^{\infty}  \int \lambda_i(\gamma) \left(\frac{d}{dt}\bigg|_{t=0} F(\gamma+tS_i)\right)^2 \, \varphi d\nu.
	\end{align*}
	Let $\xi_{\varepsilon}: \RE \rightarrow [-\varepsilon,1+\varepsilon]$ be non-decreasing such that $\xi_{\varepsilon}(t)=t$ for all $t\in [0,1]$, $0\le \xi'_{\varepsilon} \le 1$. It follows that $\xi_{\varepsilon}\circ F \rightarrow 1 \wedge F^+$. An application of the chain rule gives
	\begin{align*}
		\mathcal{E}(\xi_{\varepsilon}\circ F,\xi_{\varepsilon}\circ F) 
		&=\sum_{i=1}^{\infty} \int  \lambda_i(\gamma)\cdot \left(\frac{d}{dt}\bigg|_{t=0} \xi_{\varepsilon} \circ F(\gamma+tS_i)\right)^2 \, \varphi d\nu\\
		&= \sum_{i=1}^{\infty} \int  \lambda_i(\gamma)\cdot \xi_{\varepsilon}'(F(\gamma)) \cdot \left(\frac{d}{dt}\bigg|_{t=0} F(\gamma+tS_i)\right)^2 \, \varphi d\nu \\
		&\le \sum_{i=1}^{\infty} \int  \lambda_i(\gamma)\cdot 1\cdot\left(\frac{d}{dt}\bigg|_{t=0} F(\gamma+tS_i)\right)^2 \, \varphi d\nu \\
		&=\mathcal{E}(F,F),
	\end{align*}
	from which we derive $\mathcal{E}(1\wedge F^+,1\wedge F^+)\le \mathcal{E}(F,F)$. Thus $\mathcal{E}$ is a Dirichlet form.
	\end{proof}
\end{proposition}

\begin{proposition}
\label{local}
Let $\varphi$ satisfy \eqref{closability_cond1} and \eqref{quasi_cond}, and assume $A\ge \varepsilon \operatorname{Id}$ for some $\varepsilon>0$. Then the form $(\EF,D_Z(\EF))$ is local.
\begin{proof}
This is shown in the same manner as in Proposition 3.4 of \cite{Lobus2004} using the fact that $\nu$-a.e. implies $\varphi \nu$-a.e.
\end{proof}
\end{proposition}

\section{Quasi-regularity}
\label{section_quasiregularity}

In this section we consider weight functions $\varphi^+$ and $\varphi^-$ of the form
\begin{align}
\label{eq_def_phiplus}
	&\varphi^+(\gamma):=\exp\left\{\int_0^1\langle b_s(\gamma),d\gamma_s\rangle_{\R^d}-\hf\int_0^1|b_s(\gamma)|^2\, ds\right\},
\end{align}
and
\begin{align}
\label{eq_def_phiminus}
	&\varphi^-(\gamma):=\exp\left\{\int_0^1\langle -b_s(\gamma),d\gamma_s\rangle_{\R^d}-\hf\int_0^1|b_s(\gamma)|^2\, ds\right\},
\end{align}
where $b_s(\gamma)$ is adapted to the natural filtration of the Wiener process. Under the condition
\begin{align}
\label{eq_phi_L1}
	\frac{f}{\ \varphi^+}\in L^1(\nu), \text{ and }	\frac{f}{\ \varphi^-}\in L^1(\nu),
\end{align}
for all $f\in Z$, the closability results of Sections \ref{section_closability} and \ref{section_bilinear} hold for $\varphi^+$ as well as for $\varphi^-$.

\begin{lemma}
\label{lemma_phi_novikov}
	Let $\varphi^+$, and $\varphi^-$ be defined as in \eqref{eq_def_phiplus}, and \eqref{eq_def_phiminus}. If \eqref{eq_phi_L1} holds then we have the Novikov condition
	\begin{equation*}
		E\left[\exp{ \left\{\hf \int_0^1 |b_s(\gamma)|^2 \, ds \right\}}\right]<\infty.
	\end{equation*}
	In particular, $E[\varphi^+]=E[\varphi^-]=1$.
	\begin{proof}
		We denote 
		\begin{equation*}
			\alpha:=\exp \left\{ \int _0^1 \langle b_s(\gamma), \, d\gamma_s \rangle_{\RE^d}\right\}\text{ and }\beta:=\exp\left\{ - \hf \int_0^1 |b_s(\gamma)|^2 \, ds \right\}.
		\end{equation*}
		By \eqref{eq_phi_L1} we have $E[\alpha/\beta]<\infty$ and $E[1/(\alpha\beta)]<\infty$. Thus 
\begin{equation*}
	\infty> E\left[\left.\left(\chi_{\{\alpha<1\}}\frac{1}{\alpha}+\chi_{\{\alpha\ge 1\}}\alpha\right)\right/\beta\right]\ge E[1/\beta].
\end{equation*}
The claim follows.
	\end{proof}
\end{lemma}

\begin{proposition}
\label{quasi}
Suppose condition \eqref{eq_phi_L1} and the hypotheses of Proposition
\ref{prop_close_1}, namely conditions \eqref{closability_cond1}, \eqref{quasi_cond}, and that $A>\eps I$ for some $\eps>0$. Then, for
$\varphi=\varphi^+$ and for $\varphi=\varphi^-$, the closure of
\begin{equation*}
\EF(F,F)=\sum_{i=1}^{\infty}\int \lambda_i \left\langle S_i,DF\right\rangle_{\CM}^2\, \varphi d\nu ,\quad F\in Z,
\end{equation*}
in $L^2(\varphi \nu)$, is quasi-regular.
\begin{proof}
We prove the claim for $\varphi=\varphi^+$. The case $\varphi=\varphi^-$ is similar. To simplify notation we write $\varphi$ for $\varphi^+$. We follow \cite{Lobus2004}, see also \cite{Driver1992,Schmuland1992}.
\\\emph{Step 1}: For $r\in \N$, $l\in\{0,\ldots,2^{r-1}-1\}$, and $k=2^{r-1}+l$, set $s_k:=(2l+1)2^{-r}$. Let $x^v(p)$ denote the $v$th coordinate of $p\in  \RE^d$, $v\in \{1,\dots,d\}$. Fix $\tau \in \Omega$, $k=2^{r-1}+l$, and $v\in \{1,\dots,d\}$. Consider the functions $f_{v,k,\tau}(p):=x^v(p)-x^v(\tau(s_k)), p \in \RE^d$, and 
\begin{equation*}
F_{v,k,\tau}(\gamma):= f_{v,k,\tau}(\gamma(s_k))=x^v(\gamma(s_k))-x^v(\tau(s_k)), \quad \gamma \in \Omega;
\end{equation*}
we have $F_{v,k,\tau} \in Y$. Using the procedure of Proposition \ref{prop_close_1}(b) we get
\begin{equation}
\label{quasi_rem}
\EF(F_{v,k,\tau},F_{v,k,\tau})\le C_1 <\infty
\end{equation}
where $C_1$ the bound obtained from \eqref{ubound}.
\\\emph{Step 2}: Let 
\begin{equation*}
G_{n,\tau}(\gamma)\equiv G_{n,\tau}:=\sup_{\begin{smallmatrix} k\in\{1,\dots,n\} \\ v\in\{1,\dots,d\}\end{smallmatrix}}|F_{v,k,\tau}|, \quad n\in \N,
\end{equation*}
then
\begin{align}
\nonumber
\EF(G_{n,\tau},G_{n,\tau})
&\le \sup_{\begin{smallmatrix}
    k\in\{1,\dots,n\} \\ v\in\{1,\dots,d\}\end{smallmatrix}}
\int\sum_{i=1}^{\infty}\lambda_i(\gamma) \left\langle
  S_i,DF_{v,k,\tau}(\gamma)\right\rangle_{\CM}^2 \, \varphi d\nu\\ 
\label{eq_egtau_bound}
&\le C_1<\infty ,\quad n\in \N,
\end{align}
as in \cite{Lobus2004}, Step 2 of the proof of Proposition 3.3 together with Lemma 3.2 of \cite{Lobus2004}. We show $G_{n,\tau} \in L^2(\varphi\nu)$. Define
\begin{equation*}
	(T\gamma)_t:=\gamma_t-\int_0^t\Big(b_s^v(\gamma)\Big)_{v=1,\ldots,d} \, ds,
\end{equation*}
where $b$ is from \eqref{eq_def_phiplus} such that \eqref{eq_phi_L1} is satisfied and $b^v$ denotes the $v$th coordinate. By Lemma \ref{lemma_phi_novikov} and the Girsanov theorem $(T\gamma)_t$ is a Brownian motion under $\varphi\nu$, cf. Theorem 3.5.1 and Corollary 3.5.13 in \cite{KaratzasShreve}. We have
\begin{align*}
	&(T^{-1}\gamma)_t^v
	=\gamma^v_t+\int_0^t b_s^v(T^{-1}\gamma) \, ds
	=\left(\hf\gamma^v_t+\int_0^t \chi_{\{b_s^v(T^{-1}\gamma)>0\}} b_s^v(T^{-1}\gamma) \, ds\right)\\
	&\qquad+\left(\hf\gamma^v_t+\int_0^t \chi_{\{b_s^v(T^{-1}\gamma)\le 0\}} b_s^v(T^{-1}\gamma) \, ds\right)=:M^v_1(t)-M^v_2(t),
\end{align*}
where $M_1^v\equiv M_1^v(\gamma)$ and $M_2^v\equiv M_2^v(\gamma)$ are submartingales. For the next
calculation we remind of the fact that $\gamma$ is a random element while $\tau$ is non-random and fixed. We obtain
\begin{align*}
\notag
E\left[G_{n,\tau}^{\quad 2}\varphi\right] 
&= E\left[\left(\sup_{\begin{smallmatrix} v\in\{1,\dots,d\} \\ s\in [0,1]\end{smallmatrix}}\left|x^v\left(\gamma_s\right)-x^v\left(\tau_s\right)\right|\right)^2\varphi\right]\\
\notag
&\le 2E\left[\left(\sup_{\begin{smallmatrix} v\in\{1,\dots,d\} \\ s\in [0,1]\end{smallmatrix}}\left|x^v\left(\gamma_s\right)\right|\right)^2\varphi\right]+2\sup_{\begin{smallmatrix} v\in\{1,\dots,d\} \\ s\in [0,1]\end{smallmatrix}}\left|x^v\left(\tau_s\right)\right|^2\\
&= 2E\left[\sup_{\begin{smallmatrix} v\in\{1,\dots,d\} \\ s\in [0,1]\end{smallmatrix}} \left|x^v \left(T^{-1}\gamma_s \right) \right|^2 \right]+C_2(\tau)\\
\notag
&\le 4 E\left[\sum_{v=1}^d \sup_{ s\in [0,1]} |M^v_1(s)|^2+\sup_{ s\in [0,1]} |M^v_2(s)|^2\right]+C_2(\tau).
\end{align*}
Using Doob's inequality we get
\begin{align}
\label{eq_gncdot_doob}
E\left[G_{n,\tau}^2\varphi\right]\le 16 \sum_{v=1}^d E\left[|M^v_1(1)|^2\right]+16\sum_{v=1}^d E\left[|M^v_2(1)|^2\right]+C_2(\tau).
\end{align}
We have 
\begin{align*}
	E\left[|M^v_1(1)|^2\right]
	&=E\left[\left(\hf\gamma^v_1+\int_0^1 \chi_{\{b_s^v(T^{-1}\gamma)>0\}} b_s^v(T^{-1}\gamma) \, ds\right)^2\right]\\
	&\le \hf E[(\gamma^v_1)^2]+2E\left[\left(\int_0^1 |b_s^v(T^{-1}\gamma)| \, ds\right)^2\right]\\
	&=\hf+2E\left[\left(\int_0^1 |b_s^v(\gamma)| \, ds\right)^2 \varphi(\gamma)\right]\\
	&\le\hf+2E\left[\int_0^1 |b_s^v(\gamma)|^2 \, ds \,\varphi(\gamma)\right]\\
	&\le \hf+2E\left[\exp\left\{\int_0^1 |b_s^v(\gamma)|^2 \, ds\right\}\varphi(\gamma)\right]=\hf+2E\left[\frac{1}{\varphi^-}\right].
\end{align*}
By \eqref{eq_phi_L1} this implies
$E\left[|M^v_1(1)|^2\right]<\infty$. The same calculation holds for
\linebreak $E\left[|M^v_2(1)|^2\right]$ which gives by
\eqref{eq_gncdot_doob}
\begin{equation}
\label{eq_gncdot_bound}
E\left[ G_{n,\tau}^2\varphi\right] \le C_3(\tau),\quad n\in\N,
\end{equation}
for some positive number $C_3$ which depends on $\tau$. By \eqref{eq_egtau_bound} and \eqref{eq_gncdot_bound} we have
\begin{equation*}
\EF_1(G_{n,\tau},G_{n,\tau})< C_1+C_3(\tau)<\infty,\quad n\in \N.
\end{equation*}
\\ \emph{Step 3}: We recall again that for all $\tau\in
\Omega$ the sequence $(G_{n,\tau})_{n\in\N}$ satisfies $G_{n,\tau}\le
G_{n+1,\tau}$, $n\in\N$. By using Lemma I.2.12 of \cite{Rockner} we
observe that the function
\begin{equation*}
H_\tau(\gamma):=\sup_{\begin{smallmatrix} s\in [0,1] \\ v\in\{1,\ldots,d\}\end{smallmatrix}} \big|x^v(\gamma(s))-x^v(\tau(s))\big|, \quad \gamma\in \Omega, 
\end{equation*}
belongs to $D(\EF)$, and that $\EF(H_\tau,H_\tau)\le C_1$ as well as
\begin{equation*}
\EF_1(H_\tau,H_\tau)< C_1+C_3(\tau)<\infty.
\end{equation*}
Let $\{\tau_k:k\in \N\}$ be a dense set in $\Omega$ where we have chosen $\tau_1=\tau$. Set
\begin{equation}
\label{eq_def_kn_tau}
K_{n}:=\inf_{1\le k \le n } H_{\tau_k},\quad n\in \N.
\end{equation}
For all $n\in \N$ we have $\EF(K_n,K_n)\le C_1$ using \cite{Lobus2004}, Lemma 3.2, or \cite{Rockner}, Lemma IV.4.1. Thus
\begin{equation}
\label{eq_kn_bounded}
	\EF_1(K_n,K_n)\le C_1+C_3(\tau),\quad n\in \N.
\end{equation}
Using Lemma I.2.12 of \cite{Rockner}
we obtain $K_{n} \in D(\EF)$. Furthermore $K_{n}$, $n\in \N$, is bounded in
$(D(\EF),\EF_1^{1/2})$ by \eqref{eq_kn_bounded}. We apply the
Banach-Saks theorem in the form of  Lemma I.2.12 of \cite{Rockner}, which states that every bounded sequence in
$(D(\EF),\EF_1^{1/2})$, has a subsequence whose Cesaro means, $\bar{K}_{n_k}$, converge
strongly. Because $\|K_{n}\|_{L^2(\nu)}$ decreases in $n\in
\N$ we get
\begin{equation*}
\lim_{k\to\infty} \EF_1(\bar{K}_{n_k},\bar{K}_{n_k})=0
\end{equation*}
since $\{\tau_k\}$ is dense and $(\EF,D(\EF))$ is closed. As
$K_{n}$, $n\in \N$, is continuous, we may use Proposition 3.5  of
\cite{Rockner} Chapter III, from which it follows that there exists a
subsequence $\bar{K}_{n_l}$, $l\in\N$ of $\bar{K}_{n_k}$, $k\in \N$, and an $\EF$-nest $F_m$, $m\in \N$, such that $\bar{K}_{n_l}$ converges uniformly to zero on each $F_m$ as $l\rightarrow \infty$. Since $K_n$ is decreasing this holds even on the whole original sequence $K_n$. Let us follow an idea of \cite{Schmuland1992}, proof of Proposition 3.1. Given $\delta>0$ we can find $n$ such that $K_{n}<\delta$. We have 
\begin{equation*}
F_m\subset \bigcup_{i=1}^{n} B(\tau_i,\delta),
\end{equation*}
by \eqref{eq_def_kn_tau}, where $B(x,\delta)$ denotes the ball of radius $\delta$ centered at $x$. Now it follows that each $F_m$ is totally bounded. $F_m$ closed and totally bounded implies $F_m$ compact. Thus $F_m$, $m\in \N$ forms an $\EF$-nest consisting of compact sets.
\\\emph{Step 4}: For fixed $\tau\in \Omega$, the system of functions $F_{v,k,\tau}$, $v\in\{1,\ldots,d\}$, $k\in \N$ separates the points in $\Omega$, Step 3 showed that there is an $\EF$-nest consisting of compact sets and the form is a Dirichlet form by Proposition \ref{local}. Quasi-regularity now follows from its definition.
\end{proof}
\end{proposition}
\begin{remark}
We observe that to obtain \eqref{quasi_rem} we need the condition \eqref{quasi_cond}. We recall that by Proposition \ref{prop_close_1}(b) we have $(\EF,D_Y(\EF))=(\EF,D_Z(\EF))$.
\end{remark}
\begin{corollary}
\label{process}
There exists a diffusion process properly associated with $(\EF,D_Z(\EF))$.
\begin{proof}
The result is an immediate consequence of Propositions \ref{local} and \ref{quasi} using Theorem IV.3.5 of \cite{Rockner}.
\end{proof}
\end{corollary}

\section{Transferring the results to a geometric setting}
\label{section_geometric}
Let $M$ be a connected geometrically and stochastically complete manifold
in the sense of \cite{Hsu2002,HsuOuyang2009}. As in our main reference,
paper \cite{WangWu2008}, we assume that the manifold is torsion free. We study the form 
\begin{align}
\begin{split}
\label{geo_first}
\hat{\EF}(\hF,\hat{G}) 
&= \int \left\langle \hat{D}\hF,A\hat{D}\hat{G}
\right\rangle_\CM \, \hphi d\hnu\\
&=\int \left\langle \hat{D}\hF,\sum_{i=1}^\infty \lambda_i \left\langle
    S_i,\hat{D}\hat{G}\right\rangle_\CM S_i \right\rangle_{\CM}\,
\hphi d\hnu, 
\end{split}
\end{align}
where $\hnu$ is the Wiener measure on the path space
$P_{m_0}(M):=\{\hgamma \in C([0,1]; M):\hgamma(0)=m_0 \}$, $m_0\in
M$ and $\hphi$ is a weight function we wish to specify. We mention that similar forms and operators have been considered in \cite{Capitaine1997,WangWu2008,Hsu2002}.
\medskip

More precisely let $P_0(\R^d)\equiv\Omega:=\{\gamma \in C([0,1]; \R^d):\gamma(0)=0\}$,
$I:P_0(\R^d)\rightarrow P_{m_0}(M)$ be the It\^o map and $\hnu$ be the image
measure of the Wiener measure on $P_0(\R^d)$ under $I$. In order to recall the
construction of the It\^o map let $O(M)$ denote the orthonormal frame
bundle with respect to $M$, $\pi$ be the canonical projection $O(M)\to
M$ and $H_1,\ldots,H_d$ be the canonical horizontal vector
fields. Choose $r_0\in O(M)$ such that $\pi(r_0)=m_0$. We introduce $r_\gamma$
as the solution to the Stratonovich stochastic differential equation
\begin{equation}
\label{eq_geometric_r}
\begin{cases}
\partial r_\gamma(t)&\displaystyle =\ \sum_{i=1}^d H_i(r_\gamma(t))\partial x_i,\quad t\in [0,1],\\
r_\gamma(0)&=\ r_0,
\end{cases}
\end{equation}
$\gamma=(\gamma_1,\ldots,\gamma_d)\in P_0(\R^d)$. This defines a.e. a mapping $I:P_0(\R^d)\to
P_{m_0}(M)$ by $I(\gamma)(t):=\pi(r_\gamma(t)),\ \gamma\in P_0(\R^d),\ t\in [0,1]$, the It\^o map. We also denote by $\tang_x M$ as the tangent space of $M$ at the point $x\in M$.

The following paragraph provides compatability with
\cite{WangWu2008}. Let $K\in C([0,\infty))$ such that
\begin{equation}
\label{eq:geometric_ric:1}
\ric(X,X)\ge -K(r)|X|^2,\quad X\in \tang_xM,\ x\in B(m_0,r),\ r>0,
\end{equation}
where $B(m_0,r)$ is the geodesic ball at $m_0$ with radius $r$. Let $\rho$
denote the Riemannian distance on $M$ and $\rho_{m_0}(x):=\rho(x,m_0),\ x\in M$. 
We assume that there are constants $c_1,c_2,r_1>0$ such that the following
conditions hold,
\begin{equation*}
\hf\sqrt{(d-1)K(r)}\le c_1r, \quad r\ge r_1,
\end{equation*}
and
\begin{align}
\label{eq_73}
\begin{split}
|\ric(X,Y)|^p\le c_2 \exp\left\{\frac{1}{2}e^{-1-2c_1}\rho_{m_0}(x)²\right\}
\end{split}
\end{align}
for some $p\ge 2$ and $x\in M,\ X,Y\in \tang_xM,\ |X|=|Y|=1$.
Following the proof of Lemma 2.2 in \cite{WangWu2008} word for word we obtain
\begin{equation}
\label{eq_geometric_ric_l2}
	E\left[\int_0^1 \|\ric_{r_\gamma(t)}\|^p \, dt\right]<\infty
\end{equation}
where, as in \cite{WangWu2008},
$\ric_{r_\gamma(t)}:\R^d\to \R^d$ is defined by
\begin{equation*}
\langle
\ric_{r_\gamma(t)}(a),b\rangle:=\ric(r_\gamma(t)a,r_\gamma(t)b),\quad
a,b\in \R^d.
\end{equation*}
In particular this provides the well-definiteness in $L^2(\hnu)$ of the adjoint
directional derivative $\hat{D}_h^*$ below. We define
\begin{align}
\begin{split}
\hat{Z}:=&\Big\{\hF(\hgamma)=\hat{f}\left(\hgamma(s_1),\dots,\hgamma(s_k)\right), \hgamma \in P_{m_0}(M) : \\
& 0<s_1<\dots<s_k=1, \hat{f}\in C^{\infty}_b({M}^k), k\in \N\Big\},
\end{split}
\end{align}
and
\begin{align}
\begin{split}
\label{eq_yhat}
\hat{Y}:=&\Big\{\hF(\hgamma)=\hat{f}\left(\hgamma(s_1),\dots,\hgamma(s_k)\right), \hgamma \in P_{m_0}(M) : \\
&\hF \in \hat{Z}, s_1,\dots,s_k\in \left\{\textstyle{\frac{l}{2^n}}:l \in \{1,\dots,2^n\}\right\}, n\in \N \Big\}. 
\end{split}
\end{align}
For $\hF\in\hat{Y}$ and $\hat{f}$ as in \eqref{eq_yhat} we introduce
$(\nabla_{s_i}\hat{f})(\bar{\gamma}):=\nabla_{s_i}\hat{f}(\hgamma(s_1),\ldots,\hgamma(s_k))$
as the gradient of the function $\hat{f}$ relative to the $i$th variable
while keeping the other variables fixed. We note that $(\nabla_{s_i}\hat{f})(\bar{\gamma})$ is an element of the
tangent space $\tang_{\hgamma(s_i)}$ at $\hgamma(s_i)$. We now define
\begin{equation*}
\hat{D}_s\hF(\hgamma):=\sum_{i=1}^k s\wedge s_i \cdot r_{I^{-1}(\hgamma)}^{-1}(s_i)(\nabla_{s_i}\hat{f})(\bar{\gamma}), \quad  s\in [0,1], \hgamma \in P_{m_0}(M),
\end{equation*}
where $r$ is defined in \eqref{eq_geometric_r}. We also introduce
\begin{equation*}
\hat{\wD}_s\hF(\gamma):=\sum_{i=1}^k \chi_{[0,s_i]}(s) \Big(T^\gamma_{0\leftarrow s_i}(\nabla_{s_i}\hat{f})(\bar{\gamma})\Big),\quad s\in [0,1], \hgamma \in P_{m_0}(M).
\end{equation*}
In addition as in \cite{Houdre2003} we define the damped version
\begin{equation}
\label{eq:damped:1}
\tilde{\wD}_s\hF(\gamma):=\sum_{i=1}^k \chi_{[0,s_i]}(s) \Big(Q^*_{s_i,s}T^\gamma_{0\leftarrow s_i}(\nabla_{s_i}\hat{f})(\bar{\gamma})\Big),\quad s\in [0,1], \hgamma \in P_{m_0}(M),
\end{equation}
where $Q^{*}$ denotes the adjoint of $Q_{s,t}:\R^d\to \R^d$, which is in \cite{Houdre2003} defined by
\begin{equation*}
	\frac{dQ_{s,t}}{ds}=-\hf \ric_{r_\gamma(s)}Q_{s,t},\quad Q_{t,t}=\operatorname{Id}_{\tang_{m_0}},\quad 0\le t\le s,
\end{equation*}
or equivalently
\begin{equation}
\label{eq_def_Q}
	\frac{dQ_{s,t}}{dt}=\hf Q_{s,t} \ric_{r_\gamma(t)},\quad Q_{s,s}=\operatorname{Id}_{\tang_{m_0}},\quad 0\le t\le s,
\end{equation}
cf. also \cite{Fang1995, FangMalliavin1993}. Finally for $h\in
\CM$, $t\in \R$, $s\in[0,1]$, $\hgamma\in P_{m_0}$, and $\sigma$ being
the solution to the geometric flow equation (see \cite{Driver1992},\cite{HsuOuyang2009},\cite{Hsu2002})
\begin{equation*}
	\begin{cases}
		\frac{\partial}{\partial t}\sigma(t,s)=T^{\sigma^h(t,\cdot)(\hgamma)}_{s\leftarrow 0}r_0 h(s),\\
		\sigma^h(0,s)(\hgamma)=\hgamma(s),
	\end{cases}
\end{equation*}
we define the directional derivative
\begin{equation*}
\hat{D}_h\hF:=\lim_{t\rightarrow 0}\frac{\hF(\sigma^h(t))-\hF}{t},\quad \hF\in
\hat{Y}.
\end{equation*}
It is known that $\hat{D}_h\hF=\langle \hat{D}\hF,h\rangle_\CM$
$\hnu$-a.e. In addition we know from \cite{Capitaine1997} that the adjoint $\hat{D}_h^*$ of $\hat{D_h}$ is given by
\begin{equation*}
	\hat{D}^*_h=-\hat{D}_h+l_h,
\end{equation*}
where
\begin{equation*}
l_h(\hgamma)=\int_0^1 \left\langle \dot{h}_t+\hf\ric_{r_\gamma(t)}h_t,d\gamma_t\right\rangle
\end{equation*}
and $\hgamma=I(\gamma)$.
Having defined the geometric counterparts to the objects previously studied in this paper, analogous results hold in this new setting. We define $\hat{Z}_\CM$ to be the set of test functions of form $\sum_{i=1}^k \hat{\phi}_i(\hgamma)\cdot h_i(t)$ where $\hat{\phi}_i\in \hat{Z}$, $h_i\in \CM$ and $k\in \N$.

\begin{lemma}
\label{lemma_closable1_geometric}
	Let $M$ be specified as above. If for $p$ specified in \eqref{eq_73} we have
	\begin{equation}
		\label{eq_phi_L2p}
		\frac{\hat{g}}{\hphi} \in L^1(\hnu),\quad \hat{g}\in L^{\frac{p}{2}}(\hnu),
	\end{equation}
	then $(\hat{\mathscr{D}},\hat{Z})$ defined by
	\begin{equation*}
		\hat{\mathscr{D}}(\hF,\hat{G}):=\int \langle \hat{D}\hF,\hat{D}\hat{G}\rangle_\CM \, \hphi d\hnu, \quad \hF,\hat{G}\in \hat{Z},
	\end{equation*}
	is closable on $L^2(\hphi \hnu)$.

	\begin{proof}
		Let $\hu_n\in \hat{Z}$, 
		\begin{equation*}
			\hu_n \lni 0 \quad \text{in } L^2(\hphi \hnu), 
			\quad \text{and}\quad \hat{D}\hu_n \lni \hat{f} \quad \text{in } L^2(\hphi \hnu; \CM).
		\end{equation*}
		Choose an arbitrary $\hat{Z}_\CM\ni\hat{\psi} =\sum_{i=1}^k \hat{\phi}_i \cdot h_i$. We observe that under \eqref{eq_phi_L2p} $\hat{\psi}/\hphi\in L^2(\hphi\hnu;\CM)$. Below we use the relation $\gamma=I^{-1}(\hgamma)$. We have
		\begin{align*}
			\begin{split}
			&\int \left\langle \hat{D}\hu_n,\frac{\hat{\psi}}{\hphi} \right\rangle_\CM \, \hphi d\hnu			
			=\sum_{i=1}^k \int \left\langle \hat{D}\hu_n,h_i \right\rangle_\CM \hat{\phi}_i \, d\hnu
			=\sum_{i=1}^k \int (\hat{D}_{h_i}\hu_n) \hat{\phi}_i\,d\hnu\\
			&\qquad=\sum_{i=1}^k\int \hu_n \hat{D}^*_{h_i}\hat{\phi}_i\, d\hnu
			=\sum_{i=1}^k\int \hu_n \frac{\hat{D}^*_{h_i}\hat{\phi}_i}{\hphi}\, \hphi d\hnu\\
			&\qquad=\sum_{i=1}^k\int -\hu_n \frac{1}{\hphi} \left(\hat{D}_{h_i}\hat{\phi}_i-\int_0^1 \left\langle \dot{h}_i(t)+\hf \ric_{r_\gamma(t)}h_i(t),d\gamma_t\right\rangle\cdot\hat{\phi}_i\right)\hphi d\hnu.
			\end{split}
		\end{align*}
		This tends to $0$ as $\nti$ provided that
		\begin{equation}
		\label{eq_closability_geometric_expression}
			\frac{1}{\hphi} \left(\hat{D}_{h_i}\hat{\phi}_i-\int_0^1 \left\langle \dot{h}_i(t)+\hf \ric_{r_\gamma(t)}h_i(t),d\gamma_t\right\rangle\cdot\hat{\phi}_i\right) \in L^2(\hphi\hnu)
		\end{equation}
		for $i=1,\ldots,k$. The statement will follow then. To
                show \eqref{eq_closability_geometric_expression} we
                first observe that $(\hat{D}_{h_i}\hat{\phi}_i-\int_0^1
                \langle \dot{h}_i(t),d\gamma_t\rangle\hat{\phi}_i)/\hphi\in L^2(\hphi\hnu)$ by \eqref{eq_phi_L2p}. We denote 
		\begin{equation*}
			M_s\equiv M_s(\gamma):=\int_0^s \left\langle \ric_{r_\gamma(t)}h_i(t),d\gamma_t\right\rangle
		\end{equation*}
		and note that $M$ is a local martingale since
                $P(\int_0^1
                |\ric_{r_\gamma(t)}h_i(t)|^2\,dt<\infty)=1$ by
                \eqref{eq_geometric_ric_l2}. Using the
                Burkholder--Davis--Gundy inequality in the form of
                \cite{RevuzYor} Chapter IV, Theorem 4.1, we obtain
		\begin{align*}
                        \nonumber
				E[|M_1|^p]
				&\le C_p E\Big[\langle
                                M_1,M_1\rangle^{p/2}\Big]=C_pE\bigg[\bigg(\int_0^1
                                \big| \ric_{r_\gamma(t)}h_i(t)\big|^2
                                \,dt\bigg)^{p/2}\bigg]\\
                                \nonumber
&\le C_pE\bigg[\int_0^1 \big|\ric_{r_\gamma(t)}h_i(t)\big|^p\, dt\bigg]
		\end{align*}
for some positive constant $C_p$. With \eqref{eq_geometric_ric_l2} we
get $E[|M_1|^p]<\infty$ and from \eqref{eq_phi_L2p}
we obtain $M_1/\hphi\in L^2(\hphi\hnu)$. We have
verified \eqref{eq_closability_geometric_expression}. From \cite{Capitaine1997} we know that $\hat{Z}$ is dense in $L^2(\hnu)$ and thus we note that
$\{\hat{\psi}/\hphi:\hat{\psi}\in \hat{Z}_\CM\}$ is dense in $L^2(\hphi\hnu;\CM)$ cf. Lemma \ref{lemma_dense_testfunction}. The statement now follows.
	\end{proof}
\end{lemma}

\begin{proposition}
\label{prop_closability_geometric}
	Let $M$ be specified as above, let $\hphi$ satisfy \eqref{eq_phi_L2p}, and assume that $A\ge \varepsilon \operatorname{Id}$ for some $\eps>0$. Then Proposition \ref{prop_close_1} (a)-(d) holds for $\EF$, $Y$, $Z$, $\varphi$, $\nu$ replaced by $\hat{\EF}$, $\hat{Y}$, $\hat{Z}$, $\hphi$, $\hnu$ respectively and where $(\EF,D_Y(\EF))$ and $(\EF,D_Z(\EF))$ are replaced by the respective closures of $(\hat{\EF},\hat{Y})$ and $(\hat{\EF},\hat{Z})$ on $L^2(\hphi\hnu)$. In particular we mention that the counterpart to \eqref{quasi_cond} reads now as \eqref{quasi_geometric_cond}.
\hfill $\square$
\end{proposition}

We are interested in $\hphi(\hgamma)$ of the form 
\begin{align*}
	\hphi(\hgamma)
	&=\exp\left\{\int_0^1\left\langle \hat{V}(\hgamma_t),d\hgamma_t\right\rangle_{\tang_{\hgamma(t)}}-\hf\int_0^1 \left| \hat{V}(\hgamma_t)\right|^2_{\tang_{\hgamma(t)}}\, dt\right\}, \quad \hgamma\in P_{m_0}(M). 
\end{align*}
as $\hphi$ is then the Radon--Nikodym derivative of the diffusion measure corresponding to the general diffusion process on $M$ with generator $\hf \Delta_M+\hat{V}$, see \cite{Capitaine1997}, with respect to $\hnu$. We observe that 
\begin{equation}
	\label{eq_def_v}
	V:=r^{-1}\hat{V}\circ\pi\circ r \equiv r^{-1}\hat{V}\circ I,
\end{equation}
which says that $V\equiv V(\gamma)$ is an adapted $\R^d$-valued process. It turns out that 
\begin{align}
	\nonumber
	\hphi(\hgamma)&=\exp\left\{\int_0^1\left\langle \hat{V}(\hgamma_t),d\hgamma_t\right\rangle_{\tang_{\hgamma(t)}}-\hf\int_0^1 \left| \hat{V}(\hgamma_t)\right|^2_{\tang_{\hgamma(t)}}\, dt\right\}\\
	\label{eq_hphi_phi}
	&=\exp\left\{\int_0^1\left\langle (V(\gamma))_t,d\gamma_t\right\rangle_{\R^d}-\hf\int_0^1 \left|(V(\gamma))_t\right|_{\R^d}^2\, dt\right\}\\
	\nonumber
	&=:\varphi(\gamma), \quad \hgamma=\pi\circ r(\gamma)=I(\gamma), \gamma\in P_0(\R^d).
\end{align}
In other words $\varphi(\gamma)=\hphi(\hgamma)=\hphi\circ I(\gamma), \gamma\in P_0(\R^d)$. As in Section \ref{section_quasiregularity}, having the corresponding conditions to \eqref{eq_def_phiplus}, \eqref{eq_def_phiminus}, i.e.
\begin{align}
\label{eq_geometric_phi_L2p}
	\frac{\hat{f}}{\ \hphi^+}\in L^1(\hnu) \text{ and }	\frac{\hat{f}}{\ \hphi^-}\in L^1(\hnu)
\end{align}
for all $\hat{f}\in \hat{Z}$, we get $\hat{E}\left[\exp{ \left\{\hf \int_0^1 | \hat{V}(\hgamma_t)|_{\tang_{\hgamma(t)}}^2 \, dt \right\}}\right]<\infty$ where $\hat{E}$ denotes expectation taken with respect to $\hnu$. Thus we have the Novikov condition
\begin{equation}
\label{eq_geometric_novikov}
		E\left[\exp{ \left\{\hf \int_0^1 | (V(\gamma))_t|_{\R^d}^2 \, dt \right\}}\right]<\infty.
\end{equation}

In order to obtain quasi-regularity we sharpen the condition \eqref{eq:geometric_ric:1} and instead assume the Ricci curvature to be bounded from below, i.e. there exists some $c\in\R$, not necessarily non-negative, such that
\begin{equation}
\label{eq:ricci:1}
\ric(X,X)\ge c\|X\|^2\quad X\in \tang_x M, x\in M.
\end{equation}
We note that \eqref{eq:ricci:1}, the symmetry of the Ricci tensor, and the definition of the matrix valued function $(s,1]\ni t\to Q_{s,t}$, as the solution of \eqref{eq_def_Q}, implies that there exists $C\in\R$ such that
\begin{equation}
\label{eq_bound_Q}
\sup_{s,t\in [0,1]}\|Q_{s,t}\|\le C.
\end{equation}
We also sharpen \eqref{eq_geometric_ric_l2} and assume in addition
\begin{equation}
	\label{eq:ricci:2}
	E\left[\int_0^1 \|\ric_{r_\gamma(t)}\|^p \, dt\varphi\right]<\infty
\end{equation}
for some $p>2$. Since $\varphi\in L^1(\nu)$, \eqref{eq:ricci:2} is for example satisfied if $\ric$ is bounded on $M$.

\begin{proposition}
Suppose all the above specifications on $M$. Assume  the following conditions on $\hphi$ resp. $\varphi$, \eqref{eq_phi_L2p}, \eqref{eq_geometric_phi_L2p}, and \eqref{eq:ricci:2}. Let $A\ge \varepsilon \operatorname{Id}$ for some $\eps>0$, and 
\begin{equation}
\label{quasi_geometric_cond}
\sup_{c_1,c_2,\ldots \in \{0,1\}} \sum_{j=1}^d \sum_{p=0}^\infty \frac{\int\lambda_{i_{p,j}}\,\hphi d\hnu}{2^p} < \infty,
\end{equation}
where $i_{p,j}$ is as in Lemma \ref{lemma1}. Then
\begin{align}
\label{geo_quasi_form}
	\hat{\EF}(\hF,\hF)=\int \left\langle \hat{D}\hF,\sum_{i=1}^\infty \lambda_i \left\langle S_i,\hat{D}\hF\right\rangle_\CM S_i \right\rangle_{\CM}\,	\hphi d\hnu, \quad \hF\in \hat{Z}, 
\end{align}
is quasi-regular in $L^2(\hphi\hnu)$.

	\begin{proof}
		We proceed as we did for Proposition \ref{quasi}. For $r\in \N$, $l\in\{0,\ldots,2^{r-1}-1\}$, and $k=2^{r-1}+l$, set $s_k:=(2l+1)2^{-r}$. Let $x^v(p)$ denote the $v$th coordinate of $p\in  M$, embedded in $\R^N$, $v\in \{1,\dots,N\}$. We also write $x_s^v(\gamma)\equiv x^v(\gamma_s)$. Fix $\htau \in P_{m_0}(M)$, $k=2^{r-1}+l$, and $v\in \{1,\dots,N\}$. Consider the functions $\hat{f}_{v,k,\htau}(p):=x^v(p)-x^v(\htau(s_k)), p \in M$, and 
		\begin{equation*}
			\hF_{v,k,\htau}(\hgamma):= \hat{f}_{v,k,\htau}(\hgamma(s_k))=x^v(\hgamma(s_k))-x^v(\htau(s_k)), \quad \hgamma \in P_{m_0}(M);
		\end{equation*}
		we note $\hF_{v,k,\htau} \in D(\hat{\EF})$. We have 
		\begin{align*}
			|\langle S_i,\hat{D}\hF_{v,k,\htau}(\hgamma)\rangle_\CM|
			&=|\langle \nabla_{s_k}x^v(\hgamma(s_k)),r_{I^{-1}(\hgamma)}S_i(s_k)\rangle_{\tang_{\hgamma(s_k)}}|\\
			&=|\langle r_{I^{-1}(\hgamma)}^{-1} \nabla_{s_k}x^v(\hgamma(s_k)),S_i(s_k)\rangle_{\R^d}|.
		\end{align*}
		Now applying the Cauchy--Schwartz inequality and using the fact that we can isometrically embed $M$ in $\R^N$, cf. \cite{WangWu2008}, we obtain
		\begin{align*}
			|\langle S_i,\hat{D}\hF_{v,k,\htau}(\hgamma)\rangle_\CM|
			&\le \left|\nabla_{s_k}^N x^v(\hgamma(s_k))\right|\cdot \left|S_i(s_k)^v\right|
			\le \left|S_i(s_k)^v\right| ,
		\end{align*}
		where $\nabla^N$ denotes the gradient in $\R^N$. Thus by Proposition \ref{prop_closability_geometric} we get 
		\begin{align*}
			\hat{\EF}(\hF_{v,k,\htau},\hF_{v,k,\htau}) 
			&=\int \sum_{i=1}^\infty \lambda_i(\hgamma) \left\langle \hat{D}\hF_{v,k,\htau}(\hgamma), S_i\right\rangle_{\CM}^2\, \hphi d\hnu 
			\le C_1<\infty
		\end{align*}
		for some $C_1>0$ not depending on $v,k,\htau$, cf. Lemma \ref{lemma1} and \eqref{quasi_geometric_cond}. Let 
\begin{equation*}
\hat{G}_{n,\htau}=\sup_{\begin{smallmatrix} k\in\{1,\dots,n\} \\ v\in\{1,\dots,N\}\end{smallmatrix}}|\hF_{v,k,\htau}|, \quad n\in \N,
\end{equation*}
then
\begin{align}
\nonumber
\hat{\EF}(\hat{G}_{n,\htau},\hat{G}_{n,\htau})\le C_1<\infty ,\quad n\in \N,
\end{align}
as in Proposition \ref{quasi}. Our objective is to show $\hat{G}_{n,\htau} \in L^2(\hphi\hnu)$. I this way we obtain a bound on $\hat{\EF}_1(\hat{G}_{n,\htau},\hat{G}_{n,\htau})$ depending only on $\htau$. Then we apply the method of Proposition \ref{quasi}. Define
\begin{equation*}
	(T\gamma)_t:=\gamma_t-\int_0^t\Big((V^w(\gamma))_s\Big)_{w=1,\ldots,d} \, ds,
\end{equation*}
where $V$ is from \eqref{eq_def_v} such that \eqref{eq_geometric_phi_L2p} is satisfied and $V^w$ denotes the $w$th coordinate. By Lemma \ref{lemma_phi_novikov} and the Girsanov theorem $T\gamma_t$ is a Brownian motion under $\varphi\nu$, cf. Theorem 3.5.1 in \cite{KaratzasShreve} together with \eqref{eq_geometric_novikov} and Corollary 3.5.13 in \cite{KaratzasShreve}. We have
\begin{align*}
	&(T^{-1}\gamma)_t^w
	=\gamma^w_t+\int_0^t (V^w(T^{-1}\gamma))_s \, ds
	=\gamma^w_t+B^w_t,
\end{align*}
where $B$ is of bounded variation. Using Proposition 2.2 of \cite{Houdre2003} we have 
\begin{align*}
\notag
\hat{E}\left[\hat{G}_{n,\htau}^{\quad 2}\hphi\right] 
&\le \hat{E}\left[\left(\sup_{\begin{smallmatrix} v\in\{1,\dots,N\} \\ s\in [0,1]\end{smallmatrix}}\left|x_s^v\left(\hgamma\right)-x_s^v\left(\htau\right)\right|\right)^2\hphi\right]\\
%
%
\notag
&\le 2E\left[\left(\sup_{\begin{smallmatrix} v\in\{1,\dots,N\} \\ s\in [0,1]\end{smallmatrix}}\left|\int_0^1 E\left[\tilde{\wD}_t \{x_s^v(I(\gamma))\}\big|\F_t\right]\, d\gamma_t \right| \right)^2\varphi\right]+C_2(\htau),
\end{align*}
where $C_2(\htau):=2\left(\sup_{s\in[0,1]} \hat{E}[\left|x_s^v\right|]+\sup_{s\in [0,1]}|x_s^v(\htau)|\right)^2<\infty$ and as before $\gamma=I^{-1}(\hgamma)$. We have used $x^v\in L^1(\hnu)$, see e.g. \cite{WangWu2008}. We write 
\begin{equation}
\label{eq:alpha:1}
\alpha_s^v(t)(\gamma):= E\left[\tilde{\wD}_t \{x^v_s(I(\gamma))\}\big|\F_t\right]. 
\end{equation}
Taking into consideration that $\alpha_s(t)(\gamma)=0$ for $t>s$, we have
\begin{align}
\notag
\hat{E}\left[\hat{G}_{n,\htau}^{\quad 2}\hphi\right] 
&\le 2E\left[\left(\sup_{\begin{smallmatrix} v\in\{1,\dots,N\} \\ s\in [0,1]\end{smallmatrix}}\left|\int_0^s \alpha_s^v(t)(\gamma)\, d\gamma_t \right| \right)^2\varphi\right]+C_2(\htau)\\
\notag
%
\notag
&=2E\left[\left(\sup_{\begin{smallmatrix} v\in\{1,\dots,N\} \\ s\in [0,1]\end{smallmatrix}}\left|\int_0^s \alpha_s^v(t)(T^{-1}\gamma)\, d(\gamma+B)_t \right| \right)^2\right]+C_2(\htau)\\
\notag
&\le 4E\left[\left(\sup_{\begin{smallmatrix} v\in\{1,\dots,N\} \\ s\in [0,1]\end{smallmatrix}}\left|\int_0^s \alpha_s^v(t)(T^{-1}\gamma)\, d\gamma_t \right| \right)^2\right]\\
\label{eq_bound_hG}
&\qquad+4E\left[\left(\sup_{\begin{smallmatrix} v\in\{1,\dots,N\} \\ s\in [0,1]\end{smallmatrix}}\left|\int_0^s \alpha_s^v(t)(T^{-1}\gamma)\, dB_t \right| \right)^2\right]
+C_2(\htau).
%
\end{align}
By \eqref{eq:damped:1}, \eqref{eq_bound_Q} and \eqref{eq:alpha:1} we observe that the integrand in $\int_0^s \alpha_s^v(t)(T^{-1}\gamma)\, d\gamma_t$ is continuously differentiable according to the construction via \eqref{eq_def_Q}. Therefore the integral can path wise be treated as a Riemann--Stieltjes integral. We obtain
\begin{align}
	\notag
	&4E\left[\left(\sup_{\begin{smallmatrix} v\in\{1,\dots,N\} \\ s\in [0,1]\end{smallmatrix}}\left|\int_0^s \alpha_s^v(t)(T^{-1}\gamma)\, d\gamma_t \right| \right)^2\right]
	\le 4\sum_{v=1}^N E\left[\left(\sup_{ s\in [0,1]}\left|\int_0^s \alpha_s^v(t)(T^{-1}\gamma)\, d\gamma_t \right| \right)^2\right]\\
	\notag
	&\quad= 4\sum_{v=1}^N E\left[\left(\sup_{ s\in [0,1]}\left|\gamma_s \alpha_s^v(s)(T^{-1}\gamma)-\int_0^s \gamma_t\frac{d}{dt} \alpha_s^v(t)(T^{-1}\gamma)\, dt \right| \right)^2\right]\\
	\notag
	&\quad \le 4\sum_{v=1}^N E\left[2\sup_{ s\in [0,1]} \left(\gamma_s \alpha_s^v(s)(T^{-1}\gamma)\right)^2+2\sup_{s\in [0,1]} \left(\int_0^s \gamma_t\frac{d}{dt} \alpha_s^v(t)(T^{-1}\gamma)\, dt \right)^2\right]\\
	\notag
	&\quad \le 8\sum_{v=1}^N E\left[\sup_{ s\in [0,1]} \left(\gamma_s\right)^2\right]\cdot\left\|\sup_{s\in[0,1]}\alpha_s^v(s)(T^{-1}\gamma)\right\|_{L^\infty(\nu)} \\
	\notag
	&\qquad +8\sum_{v=1}^NE\left[\sup_{s\in [0,1]} \left(\int_0^s \gamma_t^2 \,dt \cdot \int_0^s \left(\frac{d}{dt} \alpha_s^v(t)(T^{-1}\gamma)\right)^2\, dt \right)\right]\\
	&\quad \le 8\sum_{v=1}^N 4C_3 +8\sum_{v=1}^N E\left[\left(\sup_{s\in [0,1]} \int_0^s \left(\frac{d}{dt} \alpha_s^v(t)(T^{-1}\gamma)\right)^2\, dt\right)^q\right]^{1/q} \cdot E\left[\left( \int_0^1 \gamma_t^2 \, dt\right)^p \right]^{1/p},
\end{align}
where $C_3:=\left\|\sup_{s\in[0,1]}\alpha_s^v(s)(T^{-1}\gamma)\right\|_{L^\infty(\nu)}<\infty$ according to \eqref{eq_bound_Q}. Using \eqref{eq_def_Q}, \eqref{eq_bound_Q}, and $|\left(T^{\cdot}_{0\leftarrow s}\nabla_s x^v\right) |_{\tang_{m_0}}\le 1$ since $M$ is isometrically embedded, we get
\begin{align}
\notag
	& E\left[\left(\sup_{s\in [0,1]} \int_0^s \left(\frac{d}{dt} \alpha_s^v(t)(T^{-1}\gamma)\right)^2\, dt\right)^q\right]^{1/q}\\
	\notag
	&\quad \le E\left[\sup_{s\in [0,1]} \int_0^s \left| E\left[\ric_{r_{T^{-1}\gamma(t)}} Q^*_{s,t}\circ \left(T^{\cdot}_{0\leftarrow s}\nabla_s x^v\right) (I\circ T^{-1}\gamma)\big|\F_t\right]\right|^{2q} \, dt\right]^{1/q}\\
	\notag
	&\quad \le C^2E\left[\sup_{s\in [0,1]} \int_0^s \left| E\left[\|\ric_{r_{T^{-1}\gamma(t)}}\|_{\tang_{m_0}} \big|\F_t\right]\right|^{2q} \, dt\right]^{1/q}\\
	\label{eq:bound:1}
	&\quad \le C^2E\left[\int_0^1 \|\ric_{r_{T^{-1}\gamma(t)}}\|_{\tang_{m_0}}^{2q} \, dt\right]^{1/q}
	=E\left[\int_0^1 \|\ric_{r_\gamma(t)}\|^{2q} \, dt\varphi\right]^{1/q}
	<\infty,
\end{align}
where, for the last line, we have taken into consideration \eqref{eq:ricci:2}.
For the second term in \eqref{eq_bound_hG} we have
\begin{align}
\notag
&E\left[\left(\sup_{\begin{smallmatrix} v\in\{1,\dots,N\} \\ s\in [0,1]\end{smallmatrix}}\left|\int_0^s \alpha_s^v(t)(T^{-1}\gamma)\, dB_t\right| \right)^2\right]\\
%
\notag
&\qquad \le\sum_{v=1}^N E\left[\sup_{s\in[0,1]}\left|\int_0^s |\langle\alpha_s^v(t)(T^{-1}\gamma), V(T^{-1}\gamma)_t\rangle_{\R^d} \, dt\right|^2\right]\\
\label{eq_hG_2}
&\qquad \le \sum_{v=1}^N E\left[\int_0^1 \sup_{s\in[0,1]} |\alpha_s^v(t)(T^{-1}\gamma)|_{\R^d}^2\, dt\right]\cdot E\left[\int_0^1 |V(T^{-1}\gamma)_t|_{\R^d}^2\, dt\right].
\end{align}
From \eqref{eq_bound_Q} we get
\begin{align}
\label{eq_hG_3}
&E\left[ \int_0^1 \sup_{s\in[0,1]} |\alpha_s^v(t)(T^{-1}(\gamma)|^2\, dt\right]<\infty,
%
%
\end{align}
and by \eqref{eq_hphi_phi} we get
\begin{align}
\notag
&E\left[\int_0^1 |V(T^{-1}\gamma)_t|_{\R^d}^2\, dt\right]
=E\left[\int_0^1 |V(\gamma)_t|_{\R^d}^2\, dt\varphi\right] \\
\label{eq_hG_4}
&\qquad\le E\left[\exp\left\{\int_0^1 |V(\gamma)_t|_{\R^d}^2\, dt\right\}\varphi\right]=E\left[\frac{1}{\varphi^{-}}\right]<\infty.
\end{align}
Plugging \eqref{eq_hG_3} and \eqref{eq_hG_4} in \eqref{eq_hG_2}, this together with \eqref{eq:bound:1} and \eqref{eq_bound_hG}, shows $\hat{G}_{n,\htau} \in L^2(\hphi\hnu)$ , the above formulated goal. Now having verified
\begin{equation*}
\hat{\EF}_1(\hat{G}_{n,\htau},\hat{G}_{n,\htau})\le C_1+\|\hat{G}_{n,\htau}\|_{L^2(\hphi\hnu)}^2<\infty,
\end{equation*}
the rest of the proof can be carried out in the same way as the proof of Proposition \ref{quasi}, Step 3-4. The claim follows.
	\end{proof}
\end{proposition}

\bibliographystyle{plain}
\bibliography{bibfile}

\begin{thebibliography}{10}

\bibitem{BouleauHirsch}
N.~Bouleau and F.~Hirsch.
\newblock {\em Dirichlet Forms and Analysis on Wiener Space}.
\newblock De Gruyter, Berlin, 1991.

\bibitem{Capitaine1997}
M.~Capitaine, E.~P. Hsu, and M.~Ledoux.
\newblock Martingale representation and a simple proof of logarithmic sobolev
  inequalities on path spaces.
\newblock {\em Electron. Comm. Probab.}, 2:71--81, 1997.

\bibitem{ChenWu2014}
X.~Chen and B.~Wu.
\newblock Functional inequality on path space over a non-compact riemannian
  manifold.
\newblock {\em J. Funct. Anal.}, 266(12):6753–--6779, 2014.

\bibitem{Driver1992}
B.~Driver and M.~R\"ockner.
\newblock Construction of diffusions on path and loop spaces of compact
  riemannian manifolds.
\newblock {\em C. R. Acad. Sci. Paris}, 315:603--608, 1992.

\bibitem{Fang1995}
S.~Fang.
\newblock Stochastic anticipative integrals on a riemannian manifold.
\newblock {\em J. Funct. Anal.}, 131(1):228–--253, 1995.

\bibitem{FangMalliavin1993}
S.~Fang and P.~Malliavin.
\newblock Stochastic analysis on the path space of a riemannian manifold.
\newblock {\em J. Funct. Anal.}, 118(1):249–--274, 1993.

\bibitem{Houdre2003}
C.~Houdr{\'e} and N.~Privault.
\newblock A concentration inequality on {R}iemannian path space.
\newblock In {\em Stochastic inequalities and applications}, volume~56 of {\em
  Progr. Probab.}, pages 15--21. Birkh\"auser, Basel, 2003.

\bibitem{Hsu2002}
E.~P. Hsu.
\newblock Quasi-invariance of the wiener measure on path spaces: noncompact
  case.
\newblock {\em J. Funct. Anal.}, 193(2):278--–290, 2002.

\bibitem{HsuOuyang2009}
E.~P. Hsu and C.~Ouyang.
\newblock Quasi-invariance of the wiener measure on the path space over a
  complete riemannian manifold.
\newblock {\em J. Funct. Anal.}, 257(5):1379--1395, 2009.

\bibitem{KaratzasShreve}
I.~Karatzas and S.~E. Shreve.
\newblock {\em Brownian Motion and Stochastic Calculus}.
\newblock Springer, Berlin, 2nd edition, 2000.

\bibitem{Lobus2004}
J.-U. L\"obus.
\newblock A class of processes on the path space over a compact riemannian
  manifold with unbounded diffusion.
\newblock {\em Trans. Amer. Math. Soc.}, 356(9):3751--3767, 2004.

\bibitem{Rockner}
Z.-M. Ma and M.~R\"ockner.
\newblock {\em Introduction to the Theory of (Non-symmetric) Dirichlet Forms}.
\newblock Springer, Berlin, 1992.

\bibitem{Nualart}
D.~Nualart.
\newblock {\em The Malliavin Calculus and Its Applications}.
\newblock AMS, Providence, RI, 2009.

\bibitem{RevuzYor}
D.~Revuz and M.~Yor.
\newblock {\em Continuous Martingales and Brownian Motion}.
\newblock Springer, Berlin, 1999.

\bibitem{Schmuland1992}
M.~R\"ockner and B.~Schmuland.
\newblock Tightness of general $c_{1,p}$ capacities on banach space.
\newblock {\em J. Funct. Anal.}, 108(1):1--12, 1992.

\bibitem{Stummer1993}
W.~Stummer.
\newblock The novikov and entropy conditions of multidimensional diffusion
  processes with singular drift.
\newblock {\em Probab. Theory Related Fields}, 97(4):515--542, 1993.

\bibitem{Wang2004}
F.-Y. Wang.
\newblock Weak poincaré inequalities on path spaces.
\newblock {\em Int. Math. Res. Not.}, 2004(2):90--108, 2004.

\bibitem{WangWu2008}
F.-Y. Wang and B.~Wu.
\newblock Quasi-regular dirichlet forms on riemannian path and loop spaces.
\newblock {\em Forum Math. Volume}, 20(6):1085--1096, 2008.

\bibitem{WangWu2009}
F.-Y. Wang and B.~Wu.
\newblock Quasi-regular dirichlet forms on free riemannian path spaces.
\newblock {\em Infin. Dimens. Anal. Quantum Probab. Relat. Top.},
  12(2):251--267, 2009.

\end{thebibliography}
\end{document}